\documentclass[]{article}

\usepackage{PRIMEarxiv}
\usepackage[utf8]{inputenc} 
\usepackage[T1]{fontenc}    
\usepackage{hyperref}       
\usepackage{url}            
\usepackage{booktabs}       
\usepackage{amsfonts}       
\usepackage{nicefrac}       
\usepackage{lipsum}
\usepackage{fancyhdr}       
\usepackage{graphicx}       
\graphicspath{{media/}}     
\usepackage{float}
\usepackage{amsthm,amsmath,amssymb}
\usepackage[numbers,sort&compress]{natbib}
\usepackage{doi}
\usepackage{xcolor}
\usepackage{multirow}     
\usepackage{latexsym}
\RequirePackage{fix-cm}
\usepackage{setspace}
\usepackage{parskip}
\usepackage{blindtext}
\hypersetup{
    colorlinks=true,
    linkcolor=blue,
    filecolor=blue,      
    urlcolor=black,
    citecolor=blue,
}
\title{Computational analysis for competition flows in arteriovenous fistulas based on non-contrast magnetic resonance imaging}

\author{
  Surabhi Rathore\thanks{This work was done at \textbf{Mathematical Science Group, Advanced Institute for Materials Research, Tohoku University,\\ 2-1-1 Katahira, Aoba-ku, Sendai 980-8577, Japan\\}}  \\ 
	Mathematics Area, mathLab,\\ 
        SISSA, International School for Advanced Studies,\\
        via Bonomea 265,\\
	Trieste 34136, Italy\\
 \\
   \And
	Hironobu Sugiyama \\
	Auckland City Hospital,\\
        Private Bag 92024,\\
        Auckland 1142, New Zealand \\
    \And
 \href{https://orcid.org/0000-0003-2043-3527}{\includegraphics[scale=0.08]{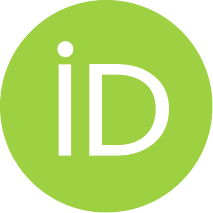}\hspace{1mm} Hiroshi Suito}\thanks{Corresponding Author. Email: \texttt{hiroshi.suito@tohoku.ac.jp}} \\
    Mathematical Science Group,\\
    Advanced Institute for Materials Research,\\ 
    Tohoku University, 2-1-1 Katahira,\\ 
    Aoba-ku, Sendai 980-8577, Japan\\
}

\begin{document}
\maketitle

\onehalfspacing

\begin{abstract}
\textbf{Introduction:}
Characteristics of hemodynamics strongly affect the patency of arteriovenous fistula (AVF) in hemodialysis patients. Because of pressure balance changes among arteries after AVF construction, regurgitating flow occurs in some patients.\\
\textbf{Methods:} 
Based on phase-contrast MRI measurements, flow types around the anastomotic site are classified to the three different types of splitting, merging, and one-way, where merging type incorporates regurgitating flow. We have performed computational simulations to analyze characteristic differences among these types. \\
\textbf{Results:}
In the merging type, a characteristic spiral flow is observed in AVF causing strong wall shear stress and large pressure drop, whereas the splitting type shows a smooth flow and gives a smaller pressure drop. The one-way case is intermediate between splitting and merging types. \\
\textbf{Conclusion:} 
Regurgitation brings about high wall shear stress near the anastomotic site because of instabilities induced by merging phenomena, for which type careful follow-up examinations are regarded as necessary.
\end{abstract}

\keywords{Arteriovenous fistula \and Hemodialysis \and Phase-contrast MRI \and Regurgitation \and Wall shear stress }
\section{Introduction}\label{sec1:intro}
End-stage renal disease (ESRD), a pathological condition in which kidney efficiency drops and low functionality persists, has been recognized worldwide as a severe chronic disease \cite{levy1996effect,hsu2008risk}. National Kidney Foundation (NKF) Kidney Disease Outcomes Quality Initiative (K/DOQI) guidelines suggest hemodialysis as the preferred treatment for ESRD patients \cite{lok2020kdoqi}. Vascular access (VA) is necessary to provide a long-patency hemodialysis. The most commonly used VA is the arteriovenous fistula (AVF) proposed by Brescia et al. \cite{brescia1966chronic}, which is created surgically between an artery and a vein. A sufficient rate of blood flow required for hemodialysis can be achieved by the pressure difference between these blood vessels. Nevertheless, reduction of the necessary blood flow rate sometimes occurs, which is responsible for AVF failure, as discussed in clinical reports \cite{van2005hemodynamics,fillinger1989beneficial}. Even with recent clinical advances, not all the precise mechanisms and factors including maturation, constriction, and occlusion are known. Although these predictive factors are debatable, most VA stenoses are known to occur near the anastomotic site of AVF. They are caused by neointimal hyperplasia \cite{stehbens1975venous}. 

The mortality rate of dialysis patients is significantly higher than that of healthy individuals, as reported by the Japanese Society for Dialysis Therapy (JSDT) \cite{mima2012hemodialysis}. Several bodies of evidence have shown AVF to be superior to other VAs such as the arteriovenous graft (AVG) and catheters in terms of patency and complications, although clinical studies \cite{allon2019vascular,viecelli2019hemodialysis} report that AVF might not necessarily be the most superior VA for elderly patients. Those with underlying diseases such as diabetes or hypertension are more likely to have blood vessel difficulties leading to primary failure of the AVF. Furthermore, primary failure in some cases might occur in AVF compared to AVG, as described in one earlier report \cite{maya2009outcomes}. To select appropriate treatments for different patient situations, algorithms for prediction that use scoring of the maturity of AVF and which use choosing between AVF or AVG have been proposed \cite{brown2017survival,voormolen2009,lok2006risk,allon2010dialysis}.

The morphology of the anastomotic sites of an AVF strongly affects the local hemodynamics within the AVF, where strong arterial blood flow dynamics exert dramatic effects on vein flow dynamics \cite{santoro2014vascular,browne2015role}. In clinical evaluations for hemodynamics through AVFs, ultrasonography is a widely available, inexpensive, and noninvasive option, although the image quality depends heavily on the skills of operators \cite{bay1998predicting,wiese2004colour,schwarz2003flow,visciano2014}. In some cases, CT scans with contrast agents, which might adversely affect chronic renal failure, are used to acquire accurate morphological information related to AVFs. Other non-invasive options are magnetic resonance imaging (MRI), which is useful without a contrast agent and which provides accurate anatomical information, as reported from an earlier study \cite{murphy2017imaging}.

Using such imaging techniques, several computational studies have been undertaken to elucidate the local hemodynamics within AVF \cite{ene2001,kharboutly2007,kharboutly2010numerical,huijbregts2007hospital,ene2012disturbed}. These studies have provided a comprehensive understanding of hemodynamics within AVFs and have revealed direct correlation between altered wall shear stress (WSS) patterns and local vessel damage. Browne et al. \cite{browne2015vivo} compared the in-vivo pressure distribution in AVF with numerical models. They reported a significant pressure drop across the fistula because of the flow instabilities at the anastomotic site. In some cases, fistula maturation depends on the anastomotic angle between an artery and a vein. This surgical factor relies on the expertise of surgeons. One earlier study \cite{stella2019assessing} has demonstrated the importance of anastomotic angles on the blood flow dynamics. That study observed more disturbed flow for the larger anastomotic angle than for the smaller anastomotic angle. Bozzetto et al. \cite{bozzetto2016transitional} investigated the flow dynamics in different AVF geometries depending on the anastomotic site location. The assumption of wall deformation would greatly complicate the problem and would require the mechanical properties of the vessel wall. However, when compliant walls are considered, the flow fields are reportedly not greatly affected compared to the rigid wall \cite{decorato2014numerical}. As described in an earlier report \cite{kikuchi2018}, some patients experience blood flow regurgitation because of changes in pressure balance between the arteries, which affects the flow dynamics at the anastomotic sites. Surabhi et al. \cite{rathore2021numerical} described detailed three-dimensional flow field structures produced because of regurgitation from another artery.

This study explores the complex flow structures, WSS, and pressure distributions using medical imaging technologies such as MRI and phase-contrast MRI (PC-MRI) with regurgitating flow from the distal artery. The stabilized finite element method is used to solve the incompressible and unsteady Navier–Stokes equations which govern blood flow. Several complex blood flow characteristics have been computed. The structure of this paper comprises the following: Section~\ref{sec2:methods} presents materials and methods, including the construction of AVF geometries using medical imaging data, computational procedures, and flow conditions. Section~\ref{sec3:results} presents the computational findings. Relevant discussion is presented in Section~\ref{sec4:Dis} with conclusions inferred from those results. 

\section{Materials and Methods}
\label{sec2:methods}
\subsection{Patient-specific data}
\label{sec2.1:PatientData}
Institutional review board (IRB) approval was obtained from Akashi Medical Center Hospital $(\#2020-12)$ and the Tohoku University School of Medicine $(\#2020-1-511)$ in accordance with the Declaration of Helsinki for the conduct of this study. We evaluated three patients who had undergone AVF surgery in their wrist area. The underlying diseases were hypertension, diabetes, nephrosclerosis, and IgA nephropathy. End-to-side anastomosis has been performed in all candidates through surgery between the end of a vein and the side of an artery \cite{konner2002anastomosis}. The evaluation was done one week after surgery. Relaxation Enhanced Angiography without Contrast and Triggering (REACT) method was used because it requires no conditions such as a high inflow rate, magnetic field homogeneity (fat suppression effect), or gated electrocardiogram. Other imaging techniques, such as Trigger Angiography Non-Contrast Enhance (TRANCE) and Time of flight angiography (TOF), which have prerequisites, are usually used for non-contrast shunt angiography, but they were not used for this study.

\begin{figure}[htbp]
\begin{center}
\includegraphics[keepaspectratio,width=8.25cm]{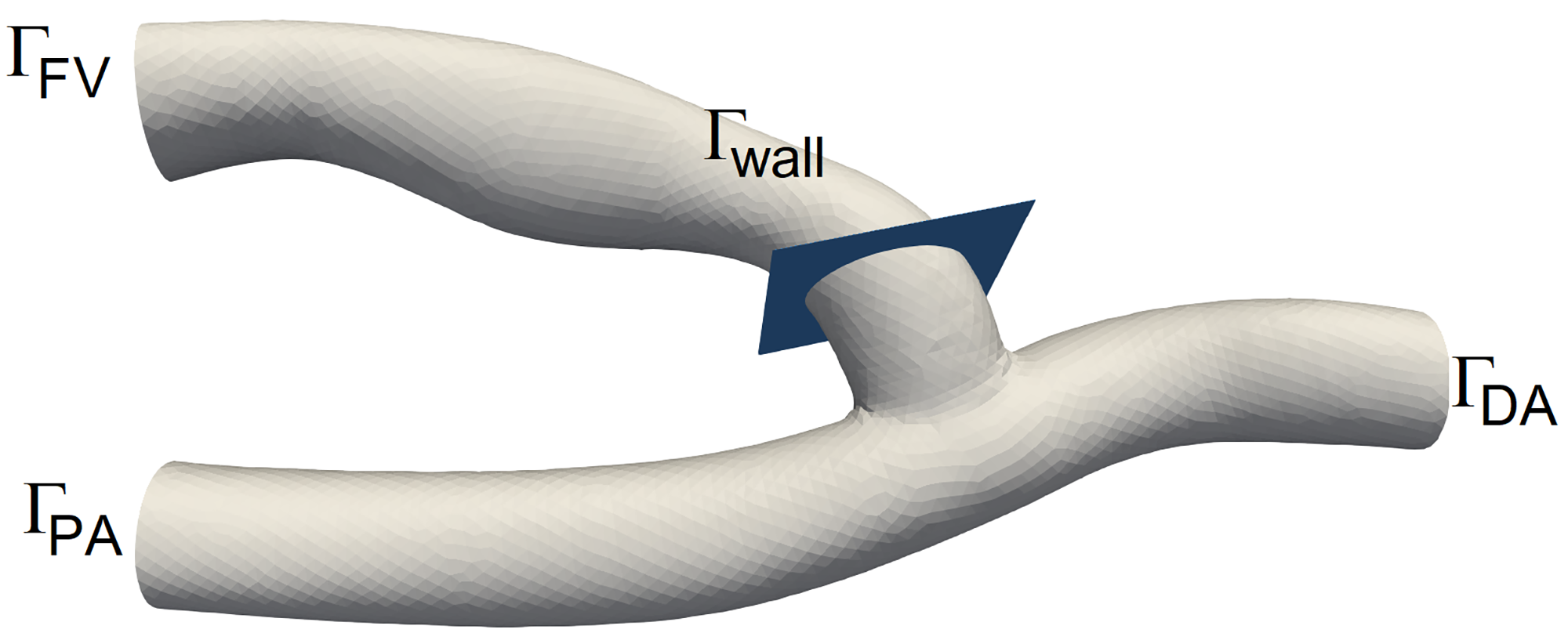}
\includegraphics[keepaspectratio,width=8.0cm]{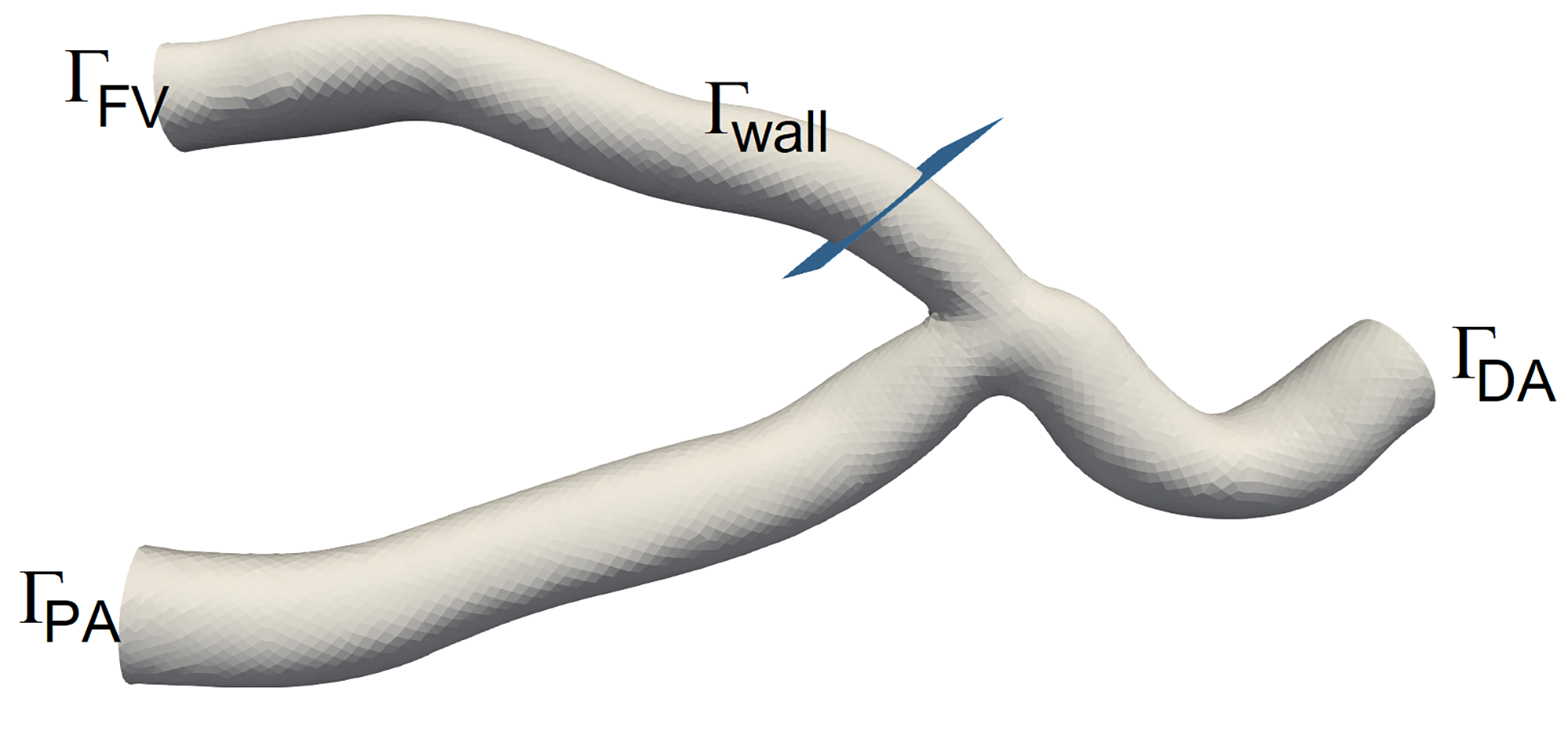}
\end{center}
\hspace{2.25cm} (a) Case S (Splitting)  \hspace{5.25cm} (b) Case M (Merging)
\begin{center}
\includegraphics[keepaspectratio,width=8.25cm]{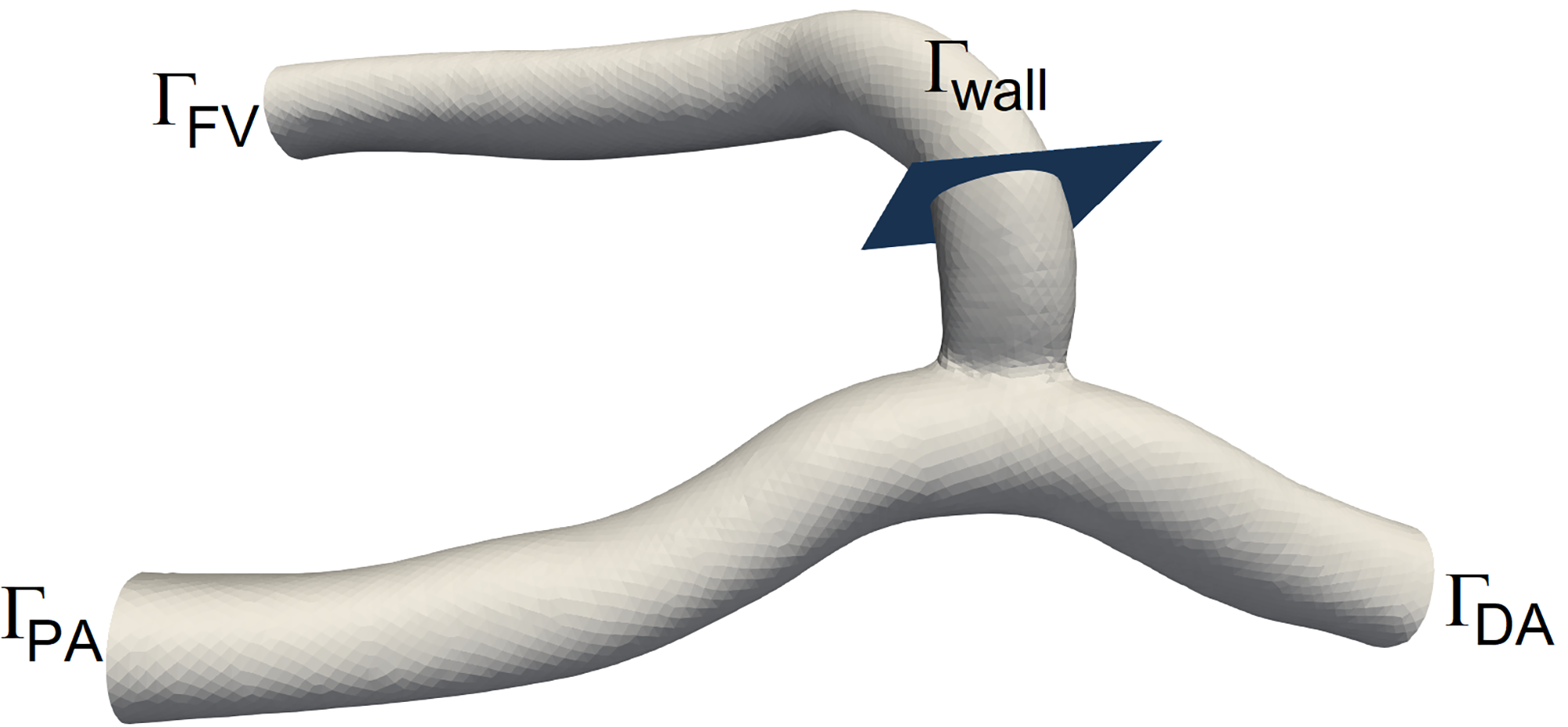}
\end{center}
\hspace{7cm} (c) Case O (One-way) 
\caption{Extracted 3D arteriovenous fistula models with planes in fistula veins.}
\label{fig1:3Dmodels}
\end{figure}

MRI medical imaging data of the targeted fistulas using the REACT method are provided in Digital Imaging and Communications in Medicine (DICOM) format. The extracted three-dimensional solid surfaces of AVF geometries are presented in Fig.~\ref{fig1:3Dmodels}. The configuration of these considered geometries consists of three parts connected at the anastomotic site presented in Fig.~\ref{fig1:3Dmodels}: the proximal artery (PA), distal artery (DA), and fistula vein (FV). Both PA and DA are parts of the radial artery. Therefore, three boundaries in our configuration other than the vessel wall boundary $\Gamma_{wall}$ are included: $\Gamma_{PA}$, $\Gamma_{DA}$, and $\Gamma_{FV}$.

\noindent We selected three patient-specific AVFs and designated these fistula geometries depending on their respective flow directions in DA as \textit{Splitting} (S), \textit{Merging} (M), and \textit{One-way} (O) cases.
\begin{itemize}
    \item \textbf{Case S:} The blood flow comes in from PA and splits into FV and DA.
    \item \textbf{Case M:} The blood flows coming in from PA and DA merge at the anastomotic site and pass out through FV. Incoming flow from DA is known as regurgitation, which occurs because of the change of pressure balance throughout the arterial system consisting of radial and ulnar arteries, and the palmar arch, by the construction of arteriovenous fistula.
    \item \textbf{Case O:} In this case, the DA flow is, quantitatively speaking, very small, indicating that the incoming flow from PA almost passes simply through FV. Aside from that, the small flow in DA causes some disturbances in the main flow from PA to FV.
\end{itemize} 

\subsection{Flow governing equations and computational procedures}
\label{sec2.2:Computation}
We have considered blood as a Newtonian fluid, the flow of which is governed by unsteady, incompressible Navier–Stokes equations given below as
\begin{align}\label{eq1:goveq}
 \rho\,\left( \mathbf{u}_{,t} + \mathbf{u} \cdot \nabla \mathbf{u} \right)  - \nabla \cdot {\mathbf{\sigma}}  =&  \mathbf{0} \hspace{2.5cm} \mbox{in} \hspace{0.25cm} \Omega\times \left(0, T\right),\nonumber
\\
\nabla \cdot \mathbf{u} =&  \mathbf{0}  \hspace{2.5cm} \mbox{in} \hspace{0.25cm} \Omega\times \left(0, T\right), \nonumber
\\
\mathbf{u} =& \mathbf{u}_{g}  \hspace{2.25cm} \mbox{on} \hspace{0.2cm} \Gamma_{g}\times \left(0, T\right), 
\\
\mathbf{\sigma} \cdot \mathbf{n}  =&  \mathbf{0}  \hspace{2.5cm} \mbox{on} \hspace{0.2cm} \Gamma_{h}\times \left(0, T\right), \nonumber
\\
\mathbf{u}\left(\mathbf{x}, 0\right) =& \mathbf{0}   \hspace{2.5cm} \mbox{on} \hspace{0.2cm}  \Omega\times \{0\}, \nonumber
\end{align}

\noindent where $\Omega \subset \mathbf{R}^{3}$ is the spatial domain with the boundary $\Gamma$ = $\Gamma_{g} \cup \Gamma_{h},$ which is defined as a disjoint union of essential and natural boundary conditions. Here, $\mathbf{u}$ and $\mathbf{n}$ respectively denote the velocity vector and the unit outward normal vector on the boundary $\Gamma$. Also, $\displaystyle \mathbf{\sigma} = \left(- p\,\mathbf{I} + 2\,\mu\,\mathbf{\epsilon} \left(\mathbf{u}\right) \right)$ is the stress tensor with the strain-rate tensor $\displaystyle \mathbf{\epsilon} \left( \mathbf{u} \right) = \frac{1}{2} \left( \nabla\,\mathbf{u} + \nabla\,\mathbf{u}^{T} \right),$ where $p$, $\mu$, and $\mathbf{I}$ respectively denote the pressure, viscosity, and identity tensor. At each boundary, flow velocity profile $\mathbf{u}_g\left(t\right)$ is given. The vessel wall is regarded as rigid.

Equation~\ref{eq1:goveq} with complex geometries is solvable numerically using Galerkin finite element methods. However, such formulations are usually associated with potential numerical instabilities \cite{donea2003finite}. Tezduyar et al. \cite{tezduyar1991stabilized,tezduyar1992incompressible} proposed stabilized finite element formulations for incompressible fluid flows. We used the streamline upwind Petrov–Galerkin (SUPG) and pressure stabilizing Petrov–Galerkin (PSPG) methods \cite{tezduyar1991stabilized,tezduyar1992incompressible} for discretizing Eq.~\ref{eq1:goveq}. Spatial domain $\Omega$ is discretized into the elements $\Omega_e = 1,2,\dots,n_{el}$, where $n_{el}$ represents the number of elements. Details of computational procedures used for this study are the same as those described for an earlier study \cite{rathore2021numerical}.

For computational simulations, blood density $\left(\rho\right)$ and blood viscosity $\left(\mu\right)$ are taken respectively as 1060.0 kg/m$^3$ and $2.66 \times 10^{-3} $ Ns/m$^2$. Periodicities of the flow fields have been examined by computing the flow fields for several heart periods to avoid effects of initial conditions because we have no clear information about them. We decided to adopt the flow fields in the third period, where the effects of the initial conditions can be regarded as negligible. Geometries of three cases have been discretized into 38,499, 493,561, and 64,455 nodes and 206,407, 263,841, and 339,797 elements, respectively, for the S, M, and O cases. Time step length of  $\Delta t = 2 \times 10^{-4}$ s  was adopted. The sparse non-symmetric system of linear equations resulting from discretization of the governing equations is solved using the GPBi-CG method, as discussed in earlier reports \cite{zhang1997gpbi,huynh2017multi}.

\subsection{Flow conditions using PC-MRI}
\label{sec2.3:BCs}
\begin{figure}[htbp]
\centering
\begin{center}
\includegraphics[keepaspectratio,width=11.5cm]{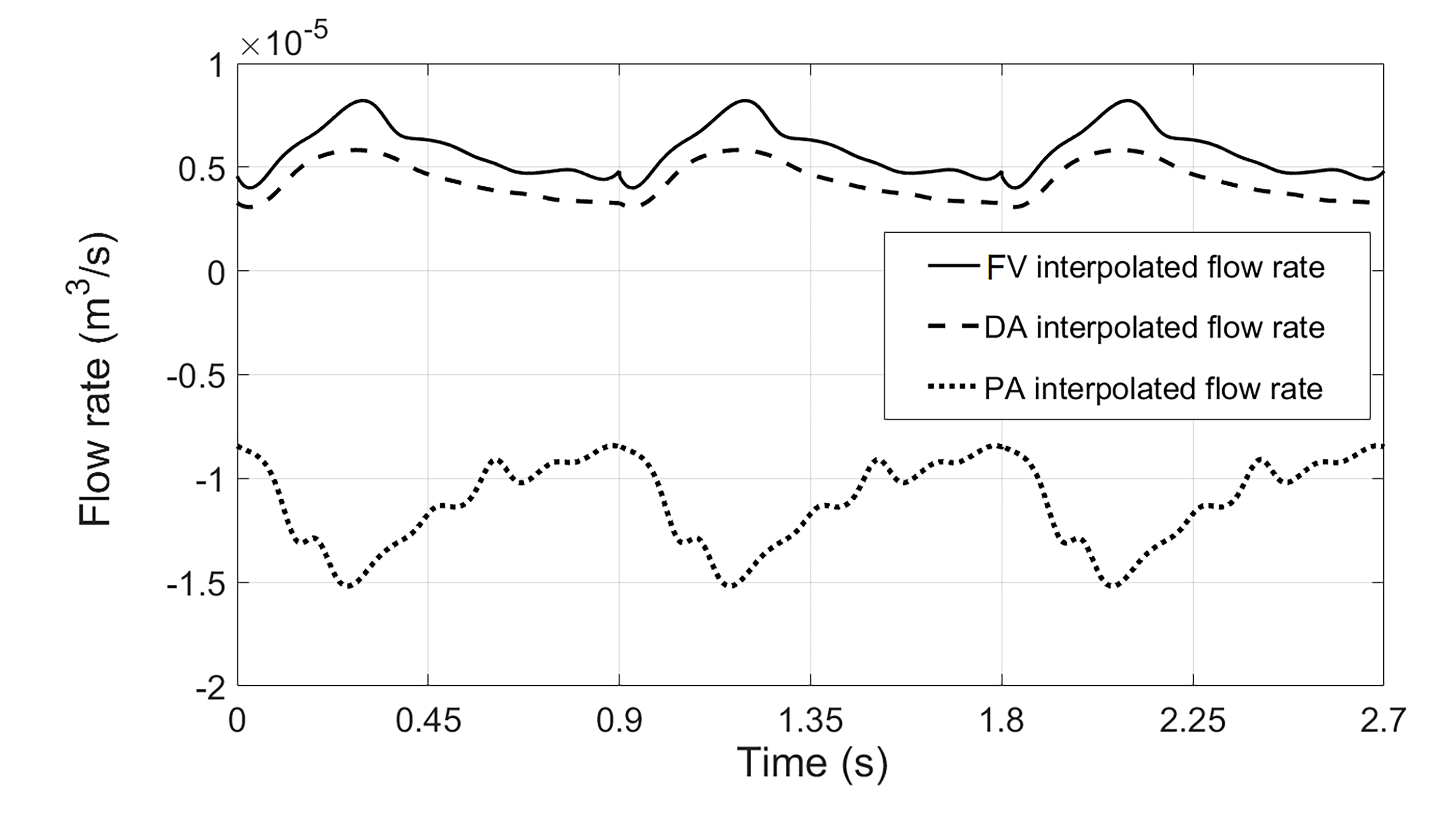}
\end{center}
(a) Case S (Splitting) 
\begin{center}
\includegraphics[keepaspectratio,width=11.5cm]{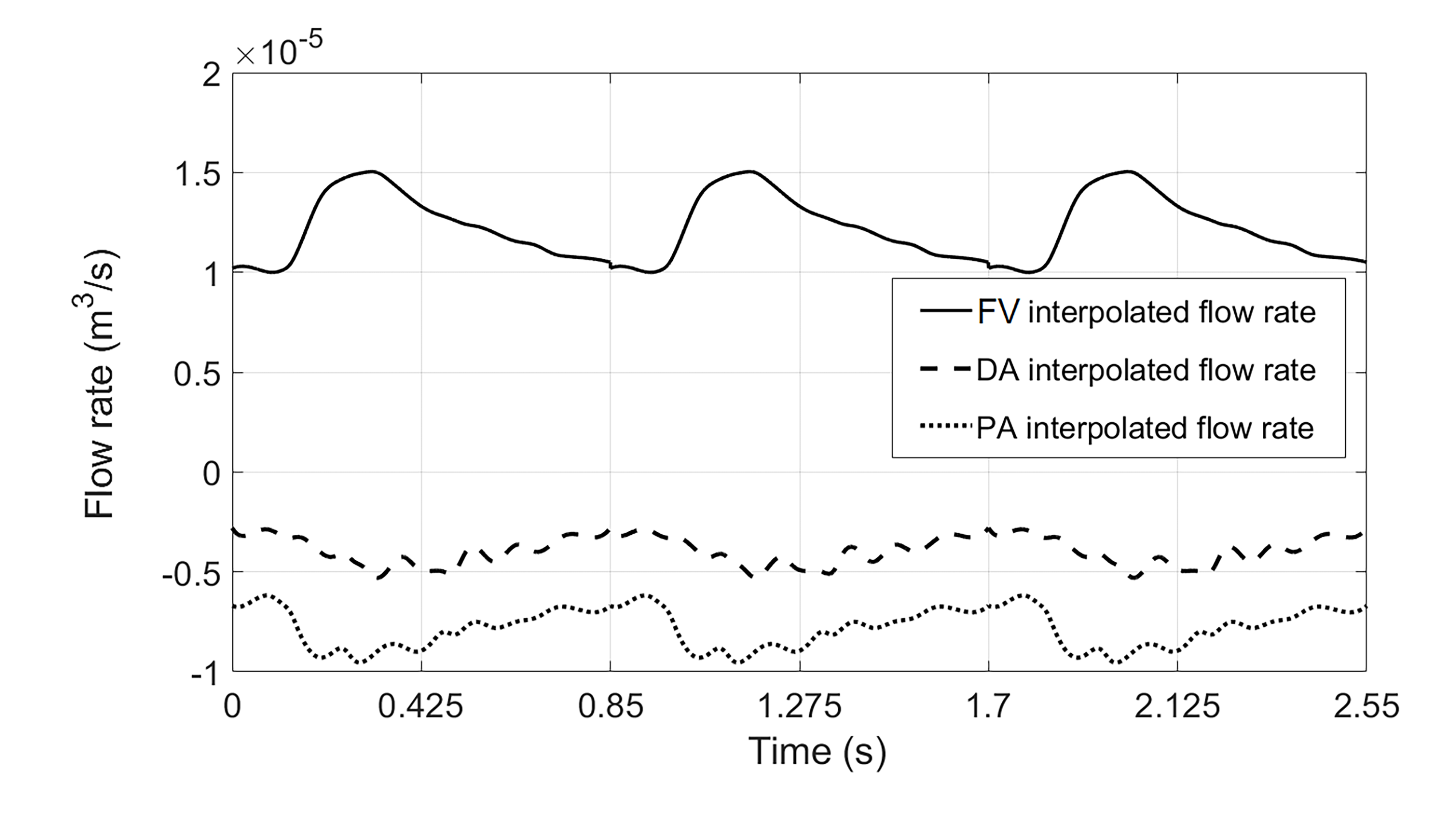}
\end{center}
 (b) Case M (Merging)
\begin{center}
\includegraphics[keepaspectratio,width=11.cm]{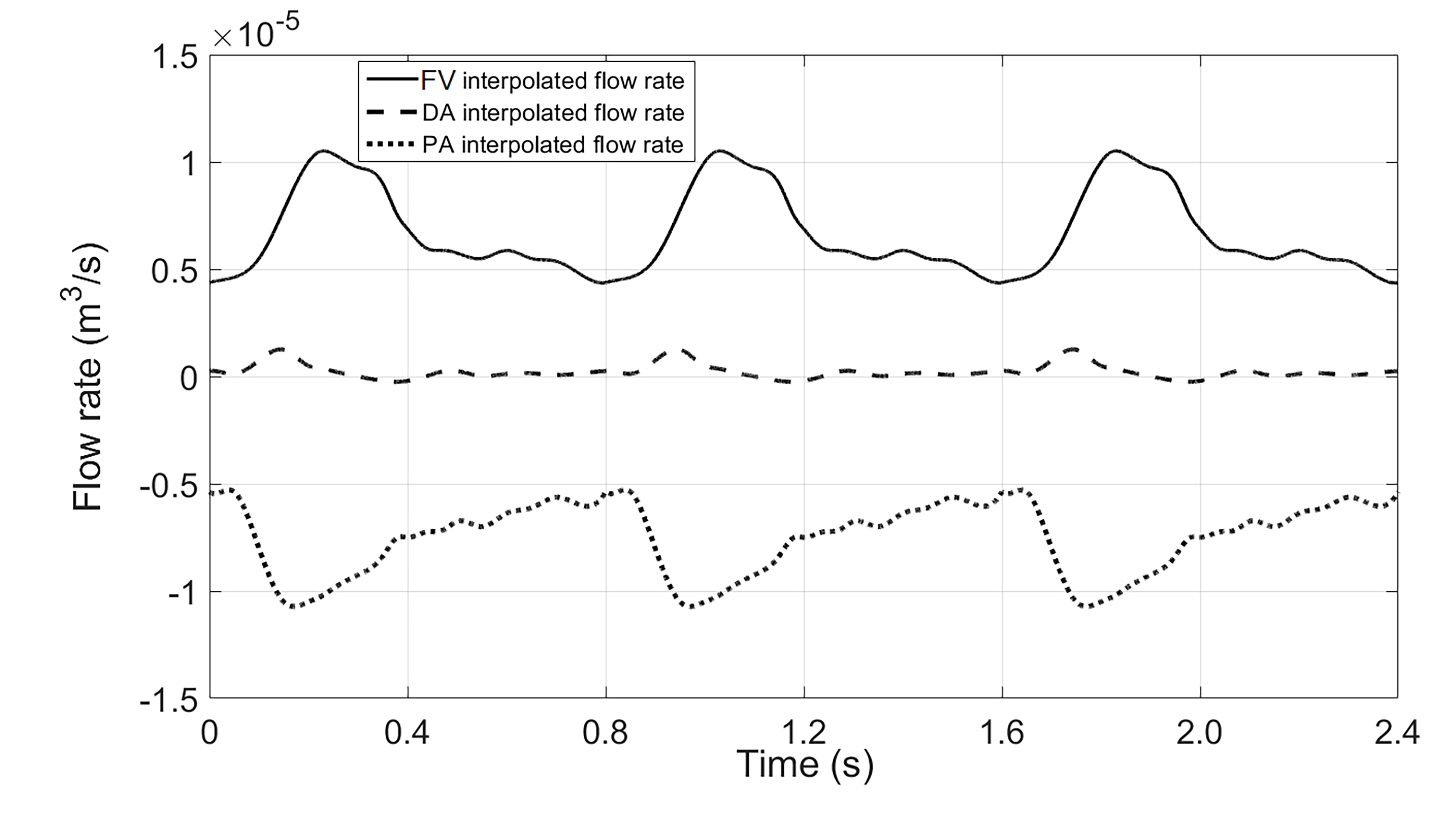}
\end{center}
 (c) Case O (One-way) 
\caption{ Flow rates for the proximal artery (PA), distal artery (DA), and fistula vein (FV) for the Splitting (S), Merging (M), and One-way (O) cases. A positive sign denotes outgoing flows; negative signifies incoming.}
\label{fig2:pcmridata}
\end{figure}

We have used PC-MRI measurement data to obtain boundary conditions for three inlet and outlet boundaries $\Gamma_{PA}$, $\Gamma_{DA}$, and $\Gamma_{FV}$. The flow rates $Q_{PA}\left(t\right)$ and $Q_{DA}\left(t\right)$ are used respectively at boundaries $\Gamma_{PA}$ and $\Gamma_{DA}$. We set $\Gamma_{FV}$ as a traction-free boundary, which means that we do not need flow rate data at $\Gamma_{FV}$. We have adopted a rectification procedure for PC-MRI measurements in this study. Because, in measurements by PC-MRI, some flow rates include certain fluctuations, which are attributable to thrilling of vessels. Such non-physiological fluctuations cause severe instabilities in computational simulations. To overcome such difficulties, we use the mass conservation law within our target vessels portion, where the (sum of) incoming flow must be equal to the (sum of) outgoing flow.
\begin{equation}\label{eq2}
    {Q_{PA}+\ Q}_{DA}+\ Q_{FV} = 0. \nonumber
\end{equation}
In fact, for M and O cases, where flow rate data at $Q_{DA}$ are fluctuating significantly, we neglect the MRI-measured flow rate profiles for DA and calculate it as
\begin{equation}
    Q_{DA}=\ {-Q}_{FV}-\ Q_{PA}, \nonumber
\end{equation}
where the negative sign of the flow rate represents incoming flow, whereas the positive sign represents outgoing flow. Rectified flow rates used as boundary conditions for computational simulations are presented in Fig.~\ref{fig2:pcmridata}.

In case S, the PA flow is divided approximately equally into FV and the DA presented in Fig.~\ref{fig2:pcmridata}(a). This flow type was discussed as Type 2 (anastomotic portion to distal) in an earlier report \cite{kikuchi2018}. In case M, flows in PA and DA traveled together to the fistula vein, as presented in Fig.~\ref{fig2:pcmridata}(b). The flow direction of DA is reversed (inward) to the anastomotic site, which might be attributable to ulnar arterial flow via the arterial arch in the palm. Such regurgitating phenomena have been discussed extensively as Type 1 (distal to anastomotic portion) \cite{kikuchi2018}. In case O, the volume of radial artery was sufficient to maintain the volume of the vein presented in Fig.~\ref{fig2:pcmridata}(c). Results show that the flow in the distal artery is very small and show that it changes the direction inlet or outlet.

\section{Results}
\label{sec3:results}
\subsection{Characteristics of flow fields in the AVFs}
\label{sec3.1:flowchar}
\begin{figure}[htbp]
\begin{center}
\includegraphics[keepaspectratio,width=8.25cm]{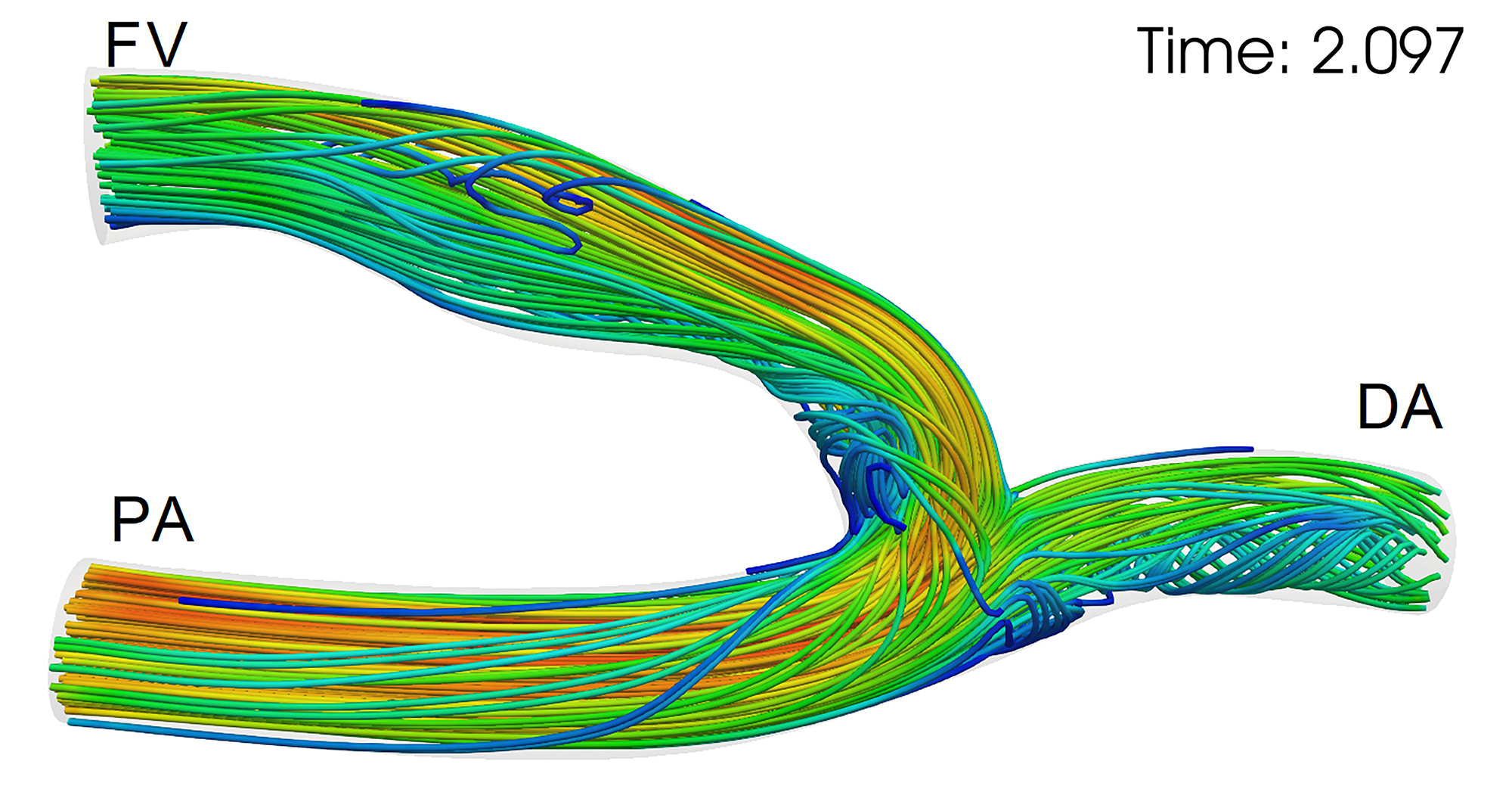}
\includegraphics[keepaspectratio,width=8.cm]{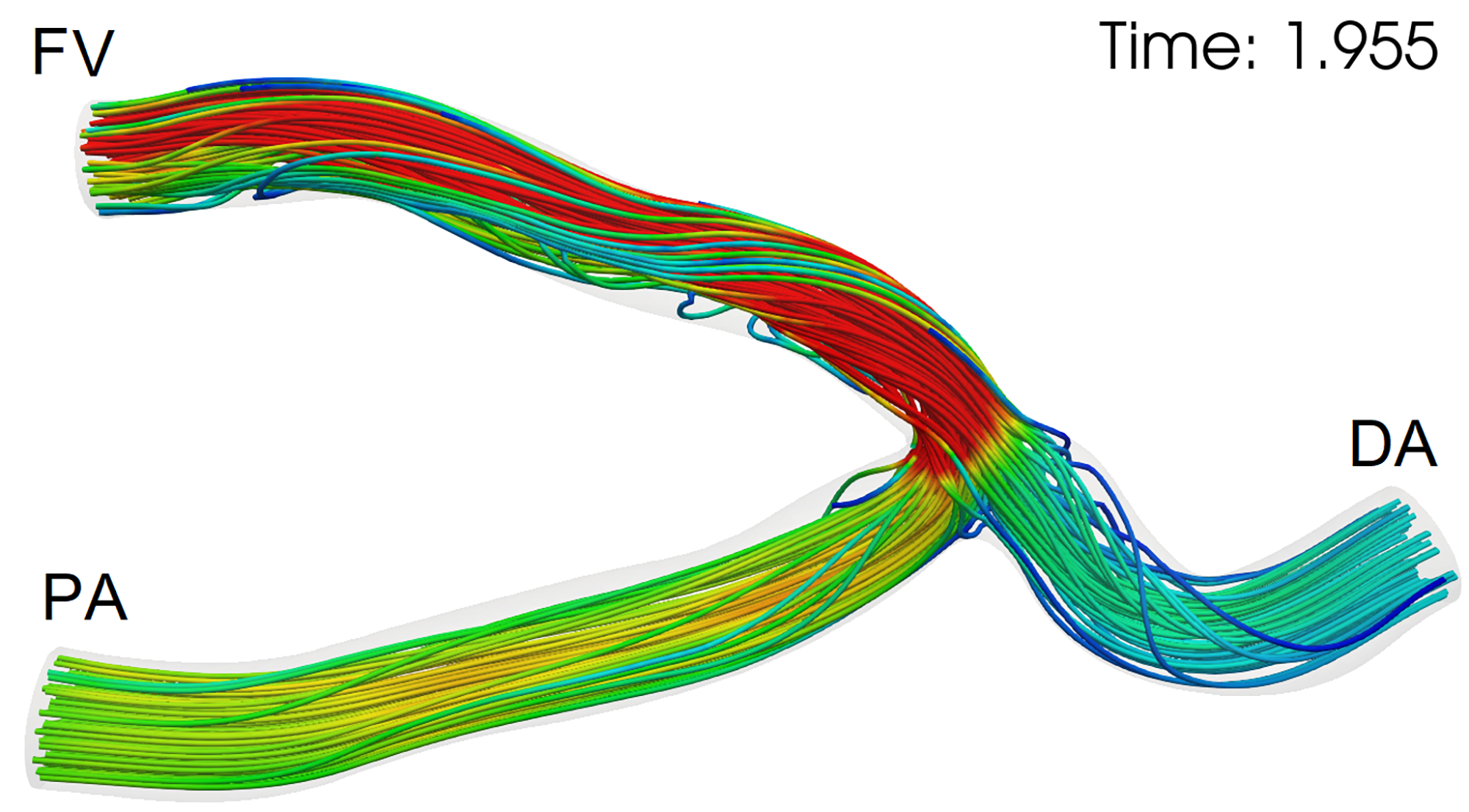}
\end{center}
\hspace{2.25cm} (a) Case S (Splitting)  \hspace{5.cm} (b) Case M (Merging)
\begin{center}
\includegraphics[keepaspectratio,width=8.cm]{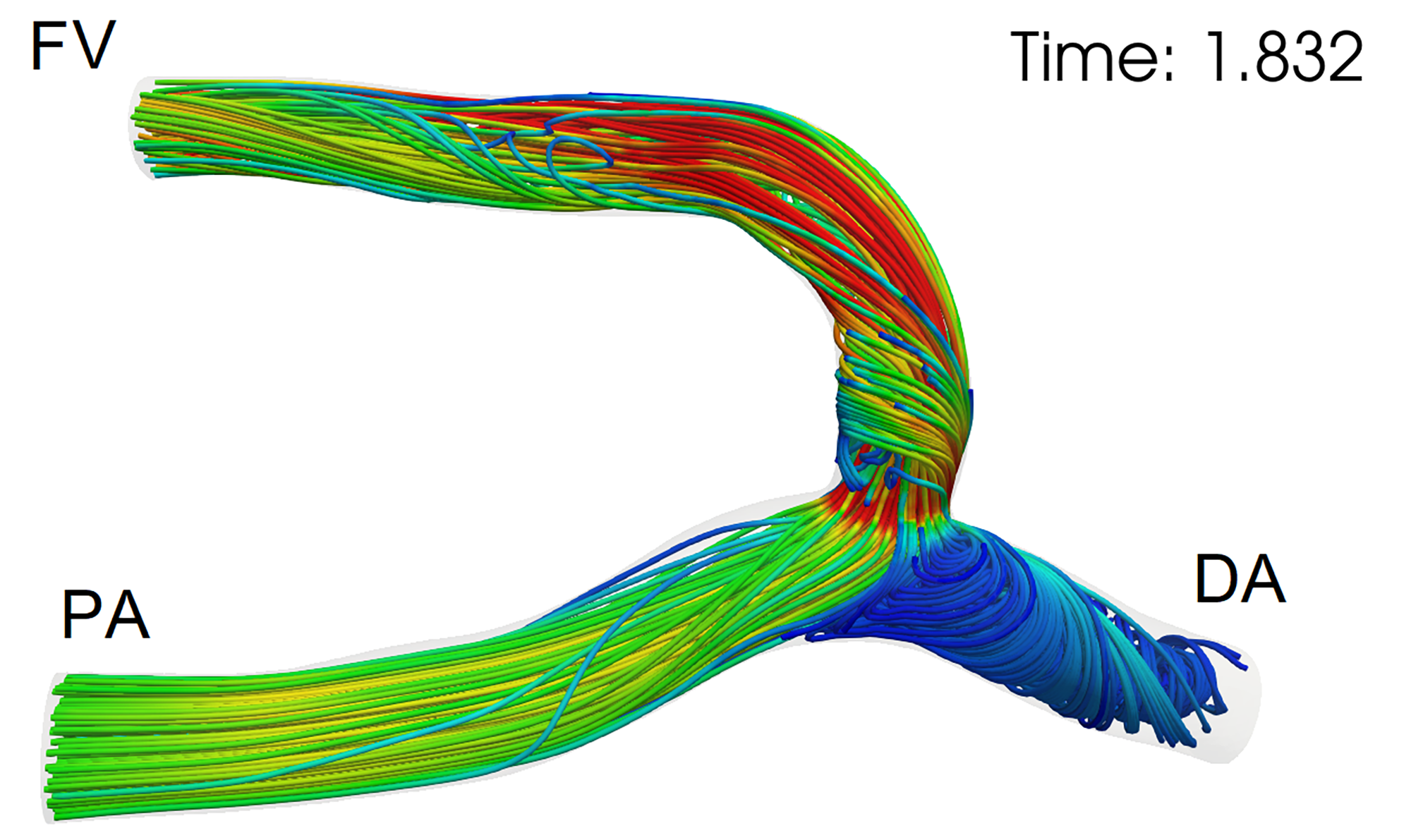}
\end{center}
\hspace{6.75cm} (c) Case O (One-way) 
\begin{center}
\includegraphics[keepaspectratio,width=7.75cm]{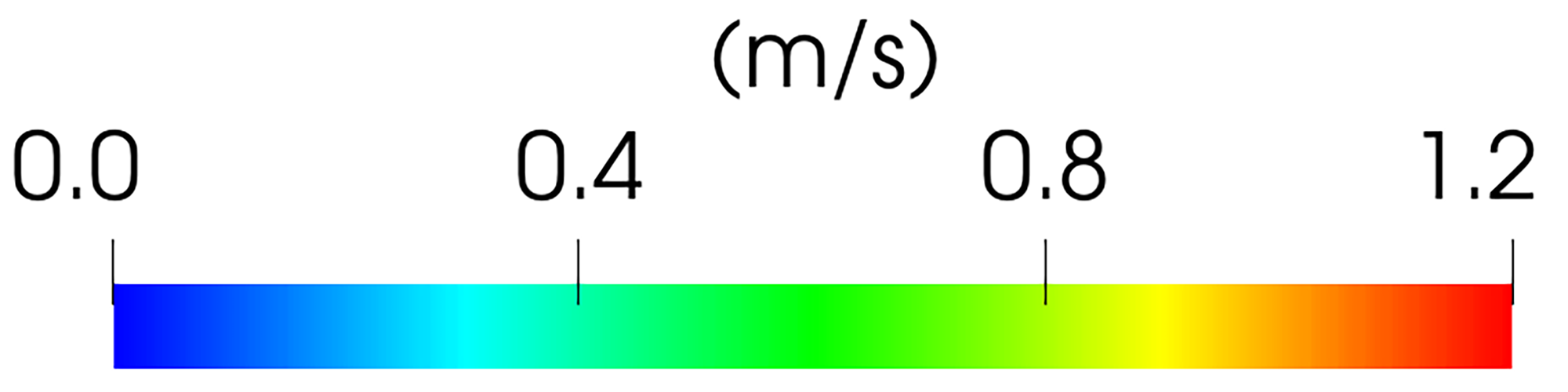}
\end{center}
\caption{Instantaneous streamlines for the S, M, and O cases.}
\label{fig3:streamlines}
\end{figure}

Figure~\ref{fig3:streamlines} shows instantaneous streamlines at the maximum inflow at $t=2.097$ s, $t=1.955$ s, and $t=1.832$ s of the third heart period, respectively, for the S, M, and O cases, where the streamline colors represent velocity magnitudes.

In case S presented in Fig.~\ref{fig3:streamlines}(a), the blood flows coming in from PA and going out from both FV and DA generate a recirculation zone at the inner side of the vein near the anastomotic site, whereas another disturbance is seen in FV because of the radius change of the vein. In case M presented in Fig.~\ref{fig3:streamlines}(b), the blood flows coming in from both PA and DA merge at the anastomotic site and go through FV. Strong flow is apparent along the FV wall in a spiral shape caused by the merging from PA and DA. In case O presented in Fig.~\ref{fig3:streamlines}(c), the blood flow is coming in mainly from PA. The small regurgitation from DA plays a role in disturbing the main flow of FV because circulation around the vessel axis is apparent in FV near the anastomotic site. We have also observed lower velocity magnitude circulation regions in DA because of marked changes in the DA flow direction.

\subsection{Wall shear stress and pressure distributions}
\label{sec3.2:wss}
\begin{figure}[htbp]
\begin{center}
\includegraphics[keepaspectratio,width=8.25cm]{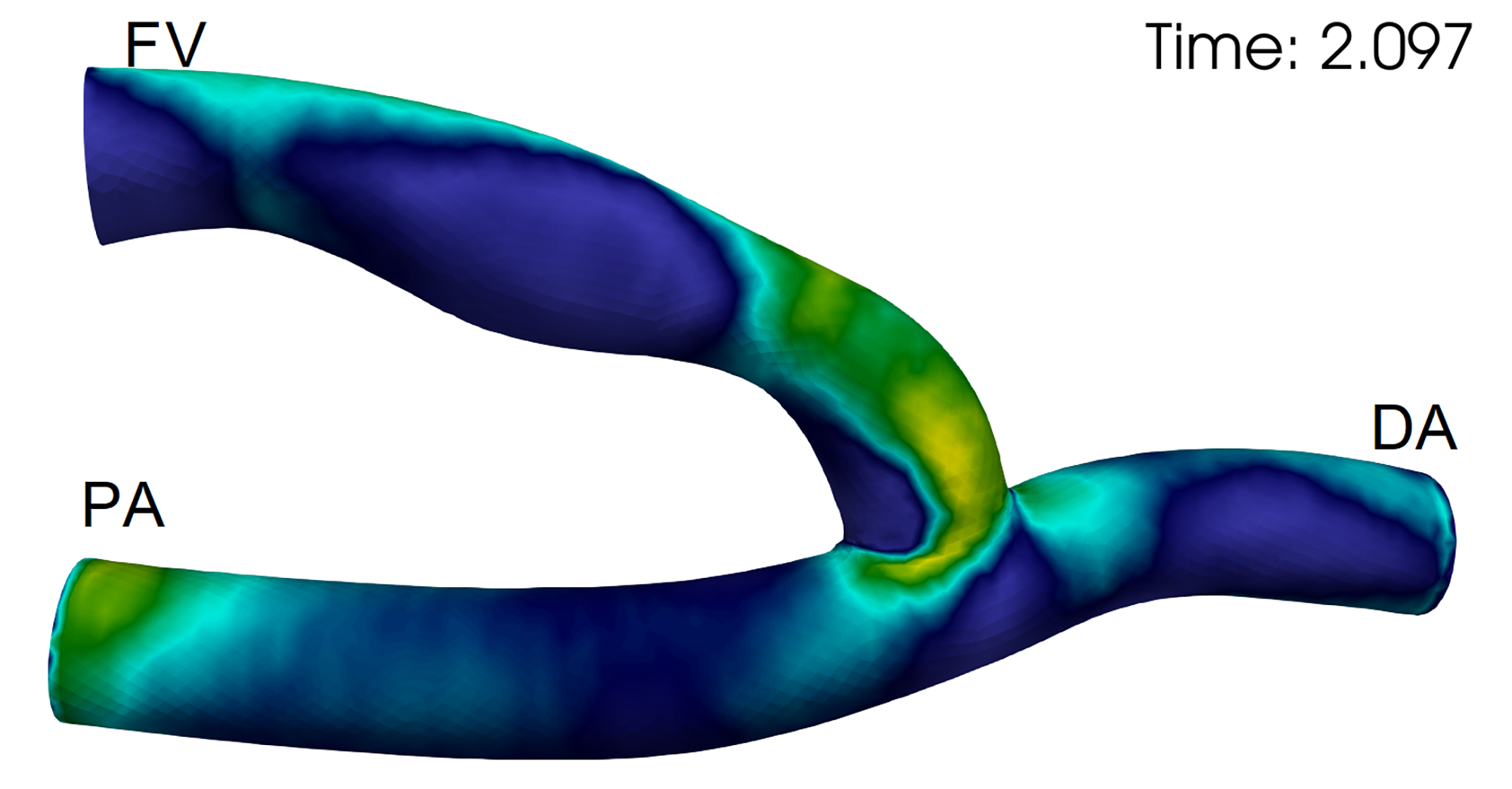}
\includegraphics[keepaspectratio,width=8.cm]{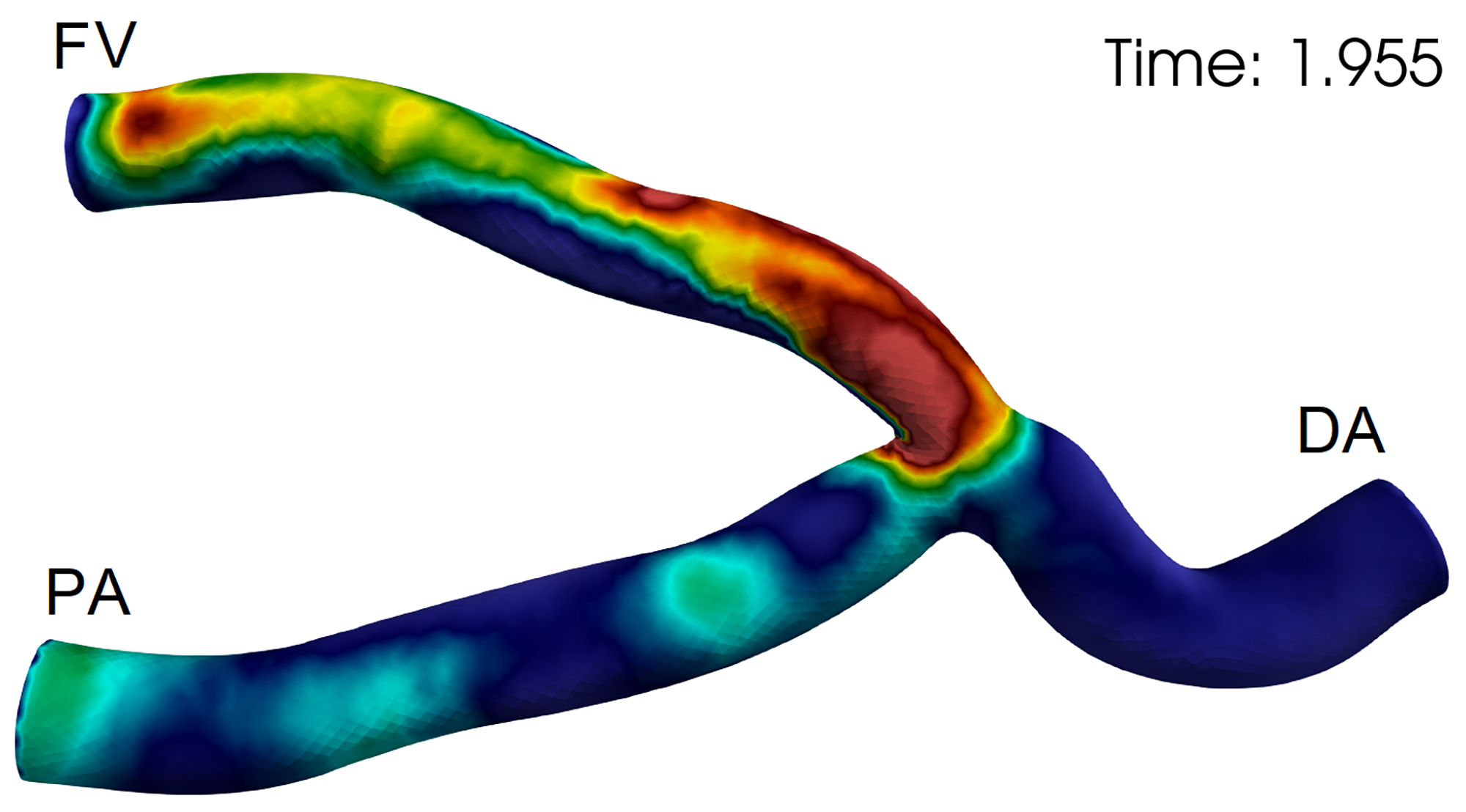}
\end{center}
\hspace{2.25cm} (a) Case S (Splitting)  \hspace{5.cm} (b) Case M (Merging)
\begin{center}
\includegraphics[keepaspectratio,width=8.cm]{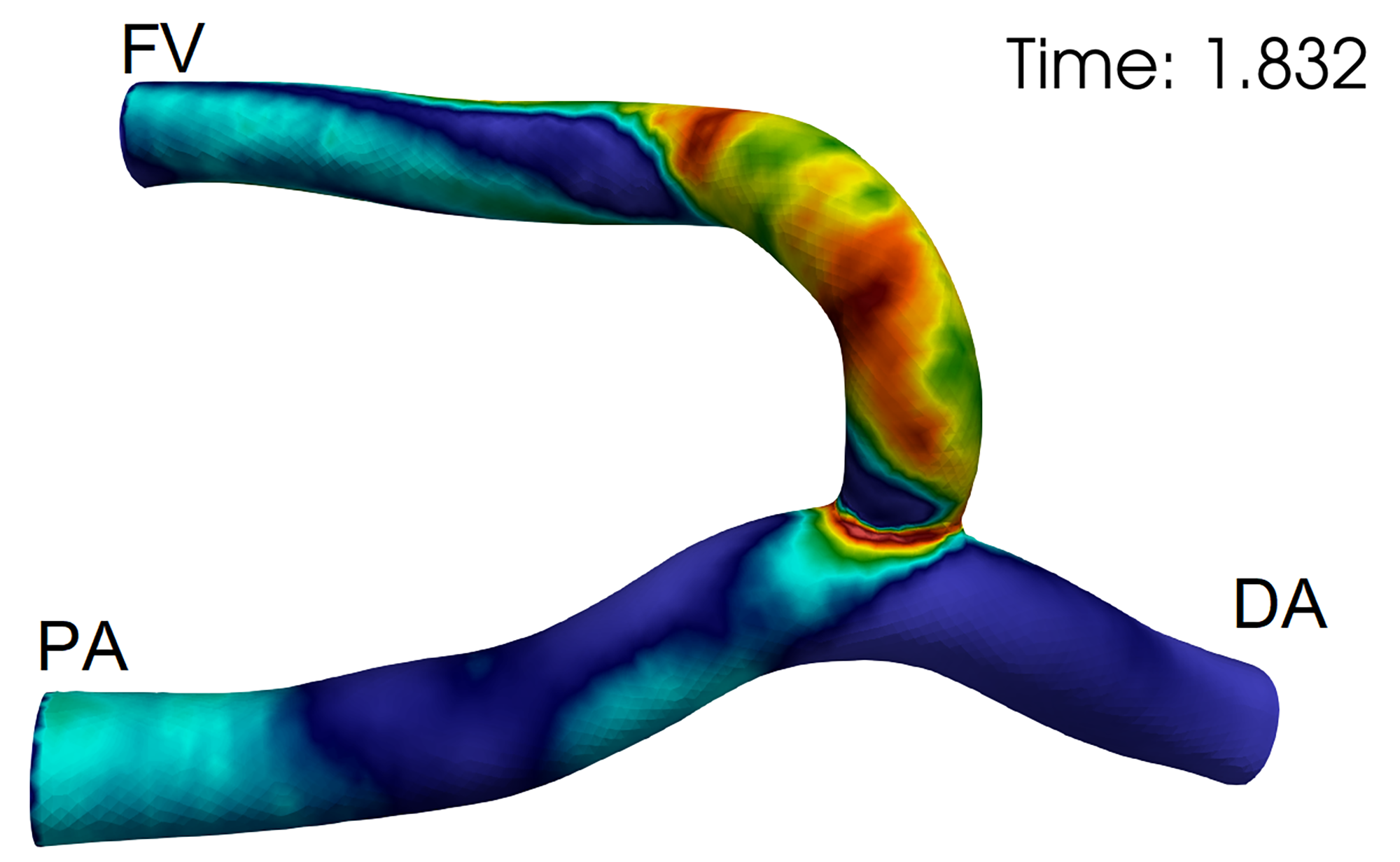}
\end{center}
\hspace{6.75cm} (c) Case O (One-way) 
\begin{center}
\includegraphics[keepaspectratio,width=7.75cm]{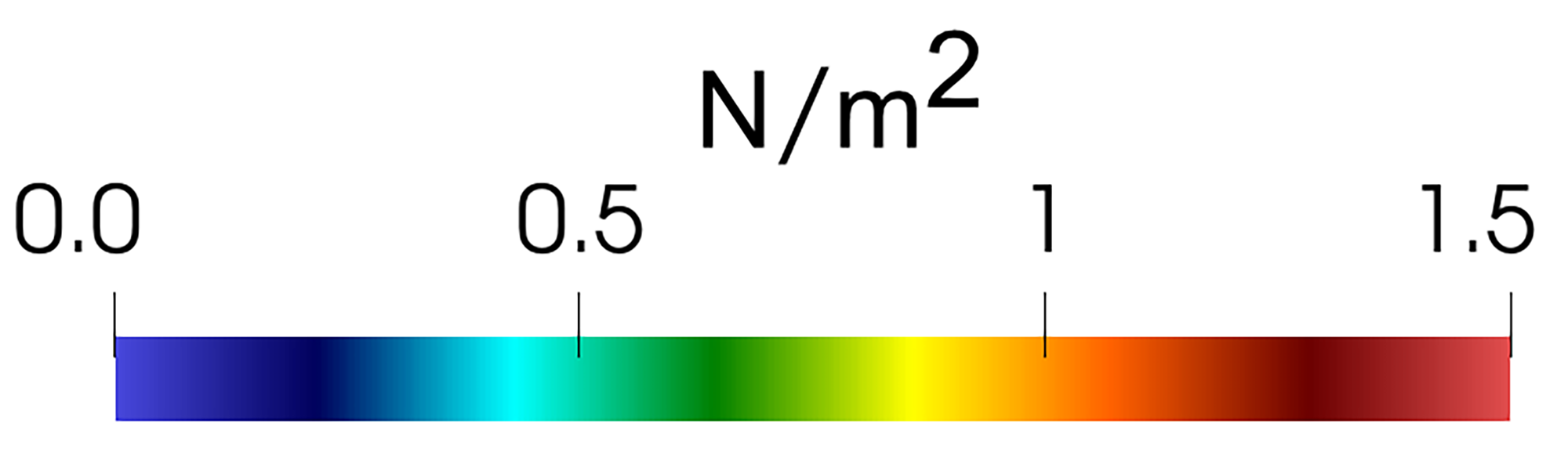}
\end{center}
\caption{Wall shear stresses for the S, M, and O cases.}
\label{fig4:wss}
\end{figure}

WSS exerted on the endothelial surface of blood vessels is defined as
\begin{equation}
  \displaystyle  \tau_{w} = \mathbf{t} - \left( \mathbf{t} \cdot \mathbf{n}\right) \mathbf{n},
\end{equation}
where $\mathbf{t}$ is a traction vector computed from the stress tensor and surface normal vector, as $\mathbf{t} = \mathbf{\sigma} \mathbf{n}$. Figure~\ref{fig4:wss} presents the WSS distributions at the time steps of the maximum inflow rates for cases S, M, and O. In case M portrayed in Fig.~\ref{fig4:wss}(b), strong WSS is apparent throughout FV which results from the strong flow merged from PA and DA. In case O portrayed in Fig.~\ref{fig4:wss}(c), strong WSS is confined near the anastomotic site which apparently results from strong circulation generated by the disturbance from DA, as depicted in Fig.~\ref{fig3:streamlines}(c). In case S portrayed in Fig.~\ref{fig4:wss}(a), WSS is maintained at a lower level, even in neighboring areas around the anastomotic site.
\begin{figure}[htbp]
\begin{center}
\includegraphics[keepaspectratio,width=5.1cm]{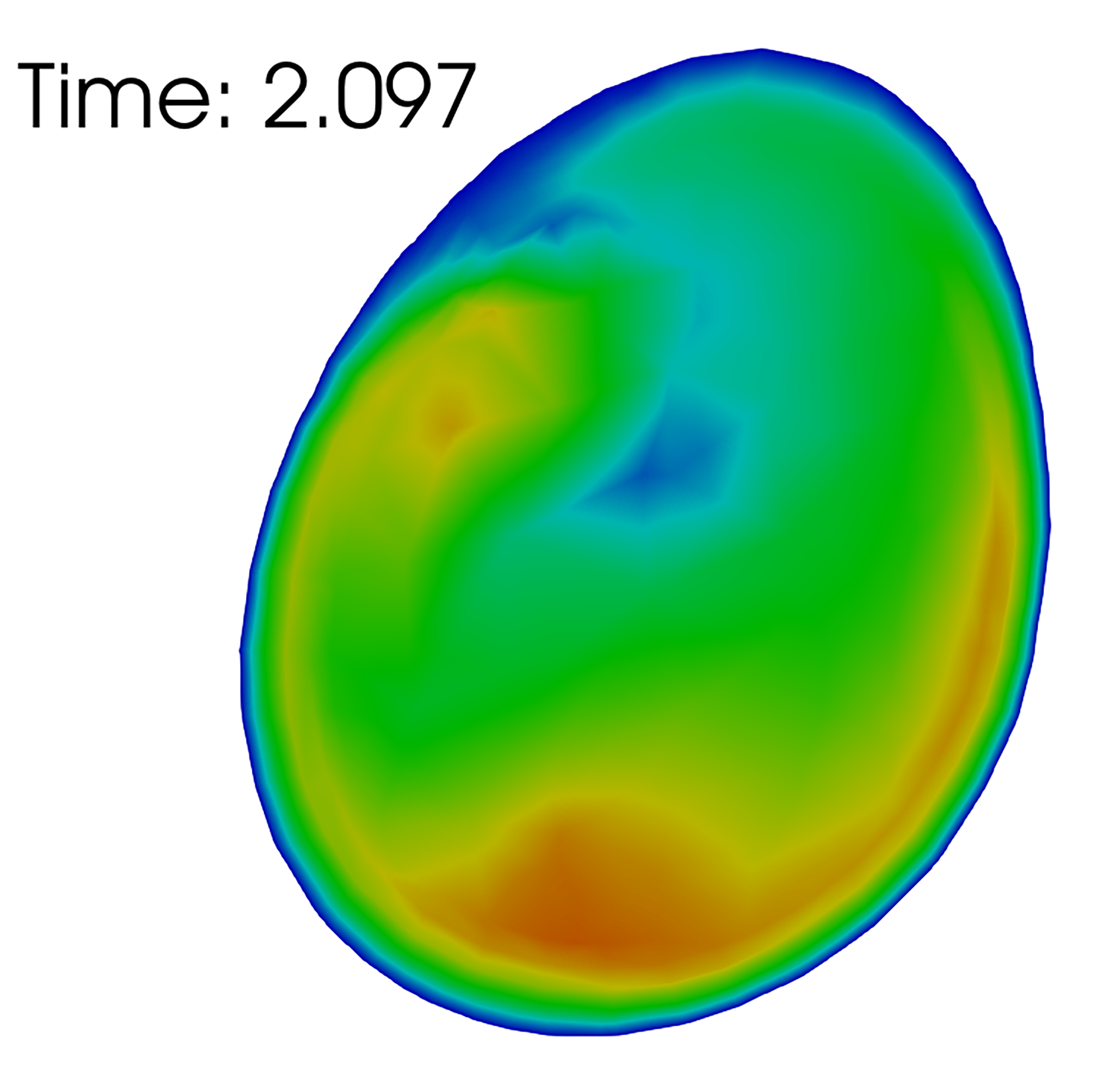}
\includegraphics[keepaspectratio,width=5.1cm]{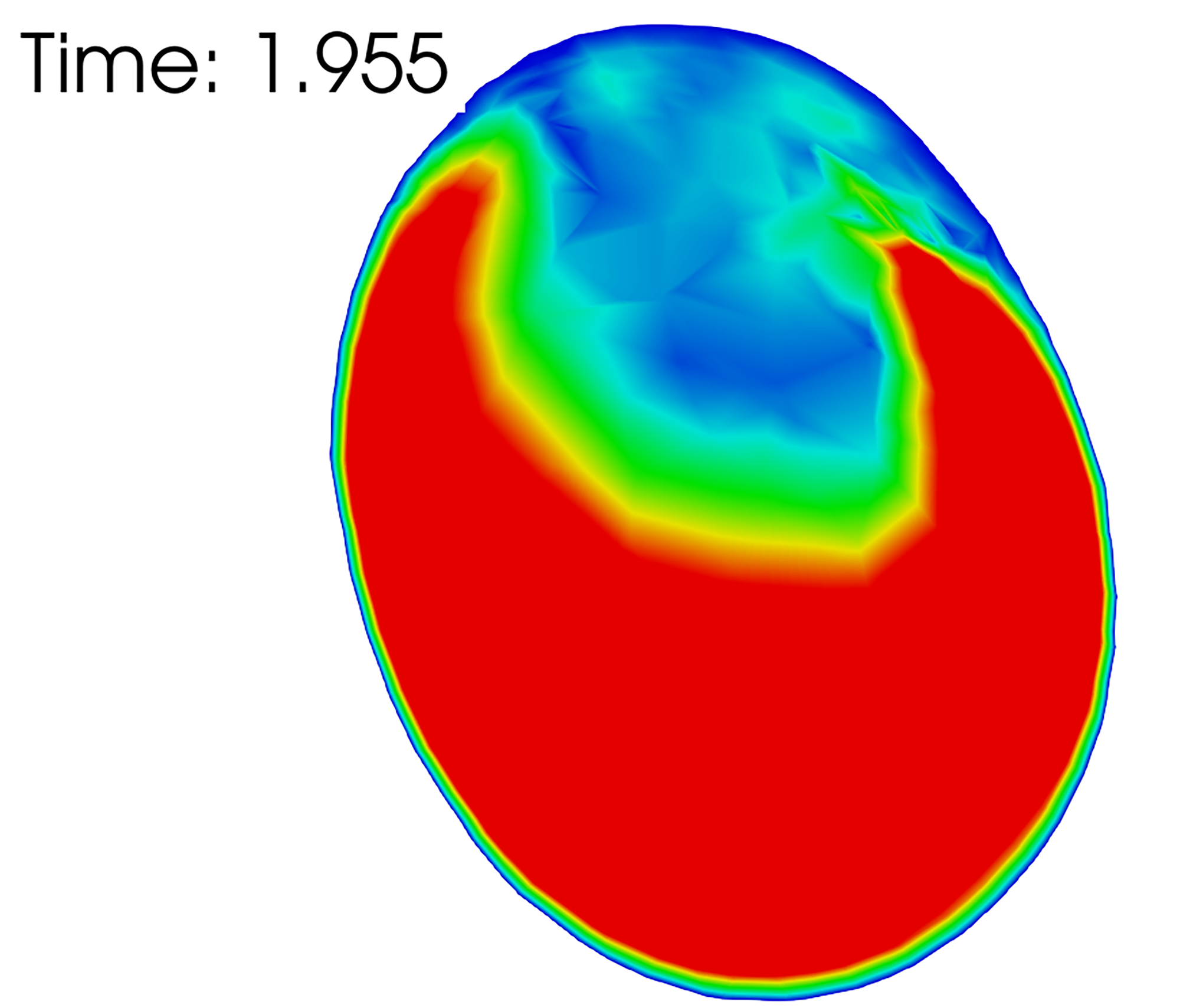}
\includegraphics[keepaspectratio,width=5.25cm]{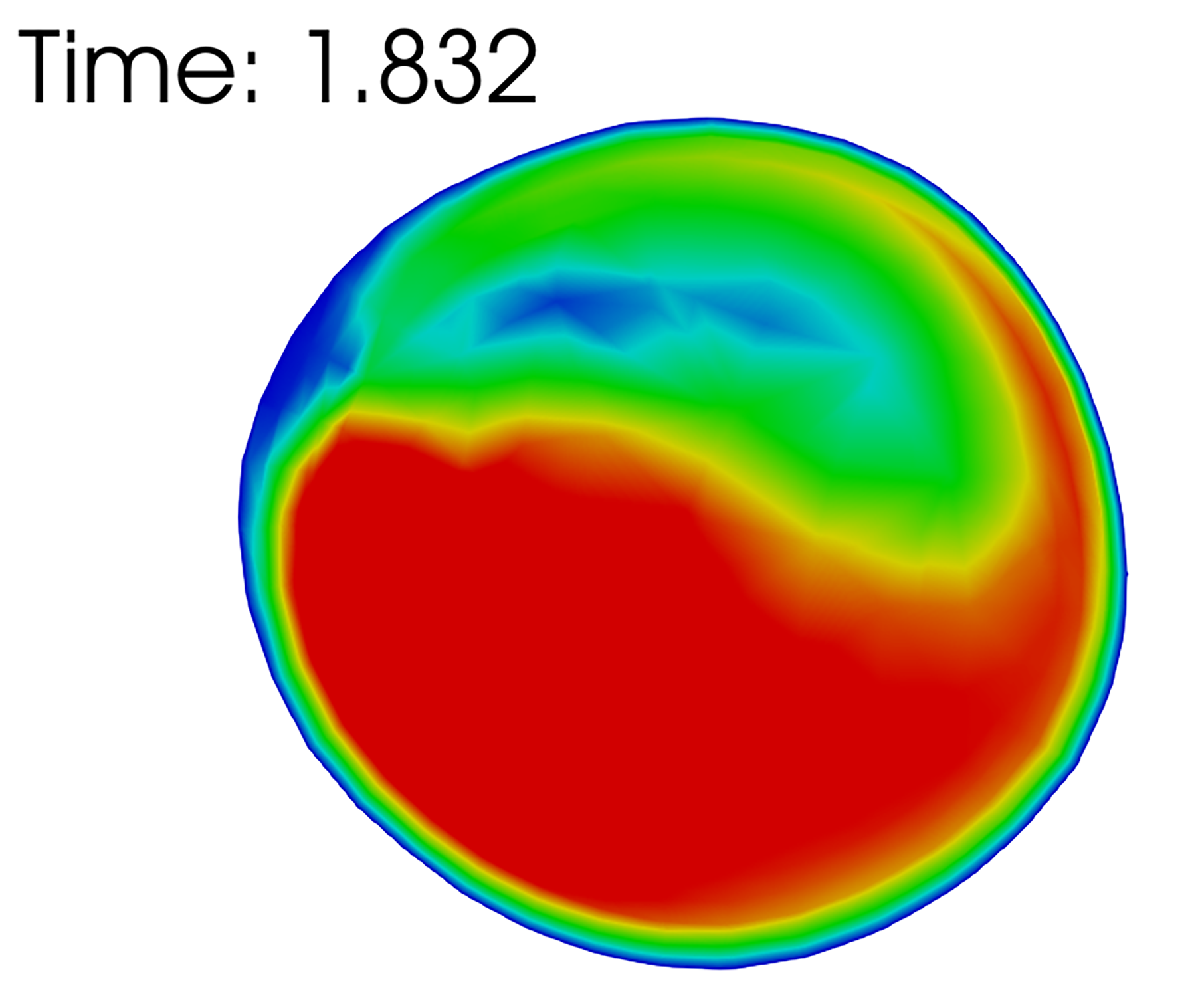}
\end{center}
\hspace{1.25cm} (a) Case S (Splitting)  \hspace{2.cm} (b) Case M (Merging) \hspace{1.75cm} (c) Case O (One-way) 
\begin{center}
\includegraphics[keepaspectratio,width=7.75cm]{Figure_3_colorbar300dpi.png}
\end{center}
\caption{Velocity magnitudes on cross sections for the S, M, and O cases.}
\label{fig5:cross}
\end{figure}

Figure~\ref{fig5:cross} presents the magnitude of velocities on the cross-sectional planes indicated in Fig.~\ref{fig1:3Dmodels} in FVs near the anastomotic sites, respectively, for cases S, M, and O. In case M in Fig.~\ref{fig5:cross}(b), a strong flow exists because of merging from PA and DA, which generates a large velocity gradient near the vessel wall and exerts strong WSS on the vessel wall. Additionally, in case O presented in Fig.~\ref{fig5:cross}(c), the larger velocity gradient is shown near the vessel wall compared with case S in Fig.~\ref{fig5:cross}(a). A critical difference between cases S and M is that the strong flow is confined in the half (lower part of Fig.~\ref{fig5:cross}(b)) of cross-section in case M, whereas the flow is distributed more uniformly over the cross-section in case S. Case O is intermediate between cases M and S.
\begin{figure}[htbp]
\begin{center}
\includegraphics[keepaspectratio,width=8.25cm]{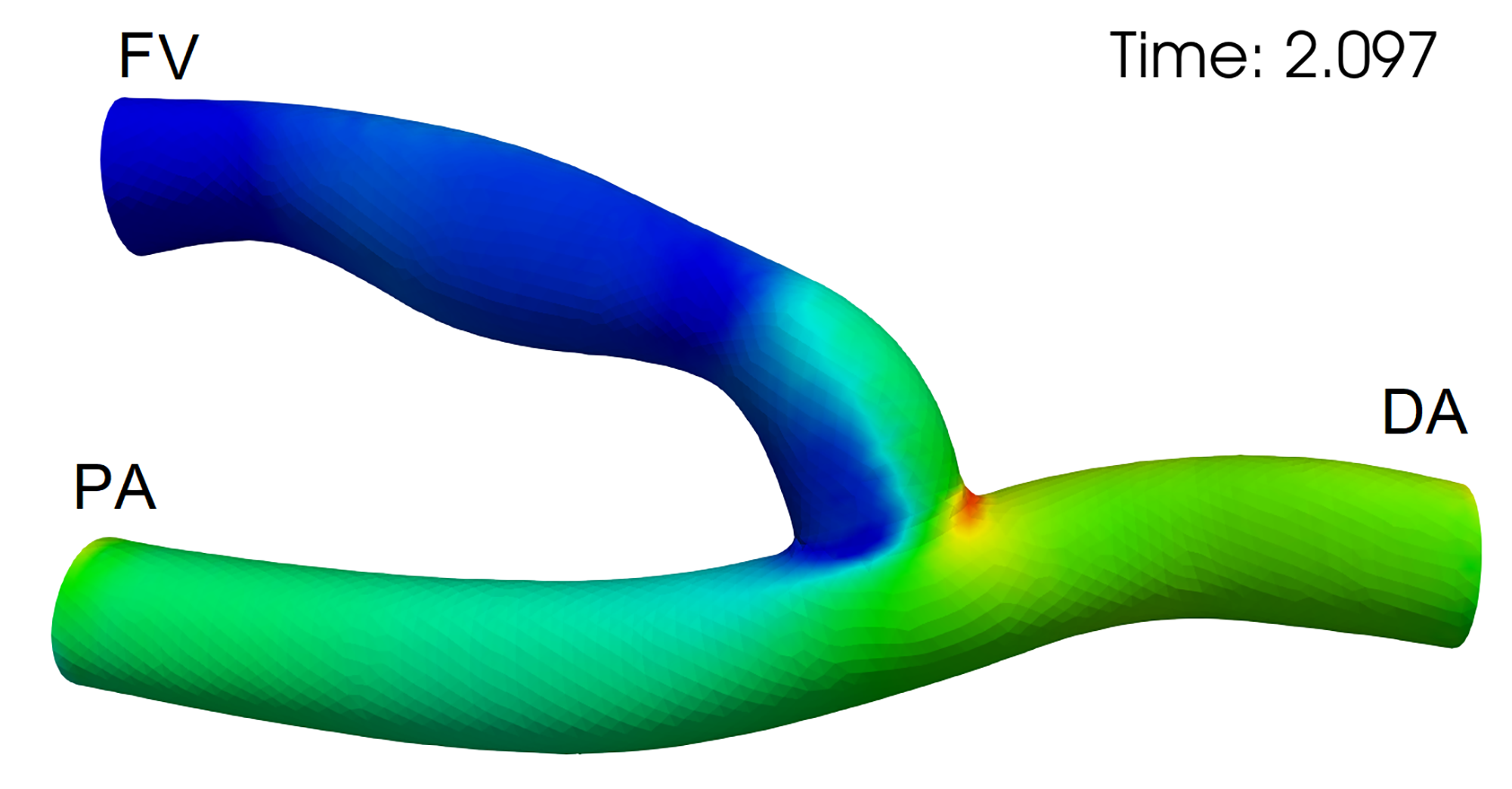}
\includegraphics[keepaspectratio,width=8.1cm]{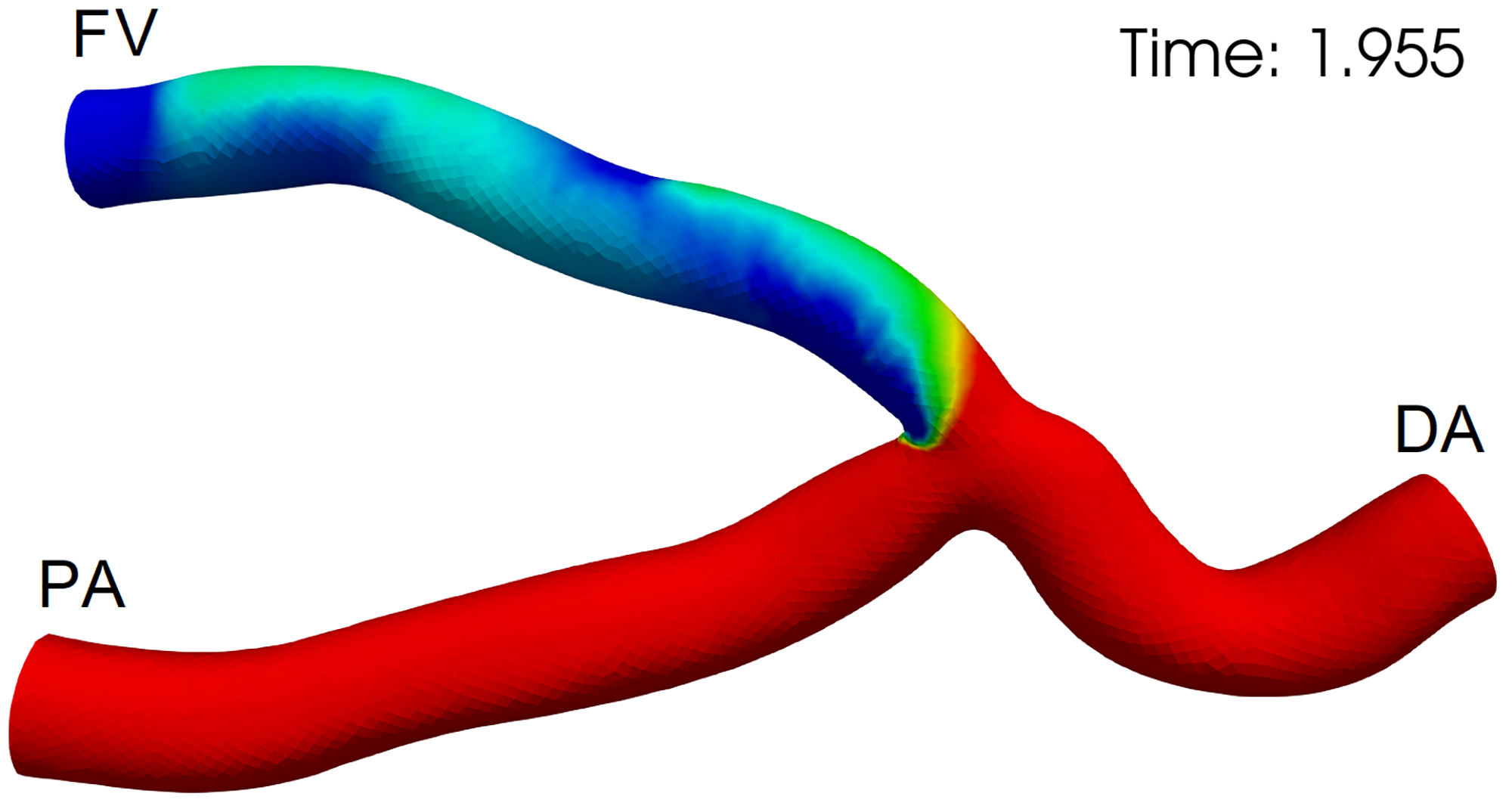}
\end{center}
\hspace{2.25cm} (a) Case S (Splitting)  \hspace{5.cm} (b) Case M (Merging)
\begin{center}
\includegraphics[keepaspectratio,width=8.cm]{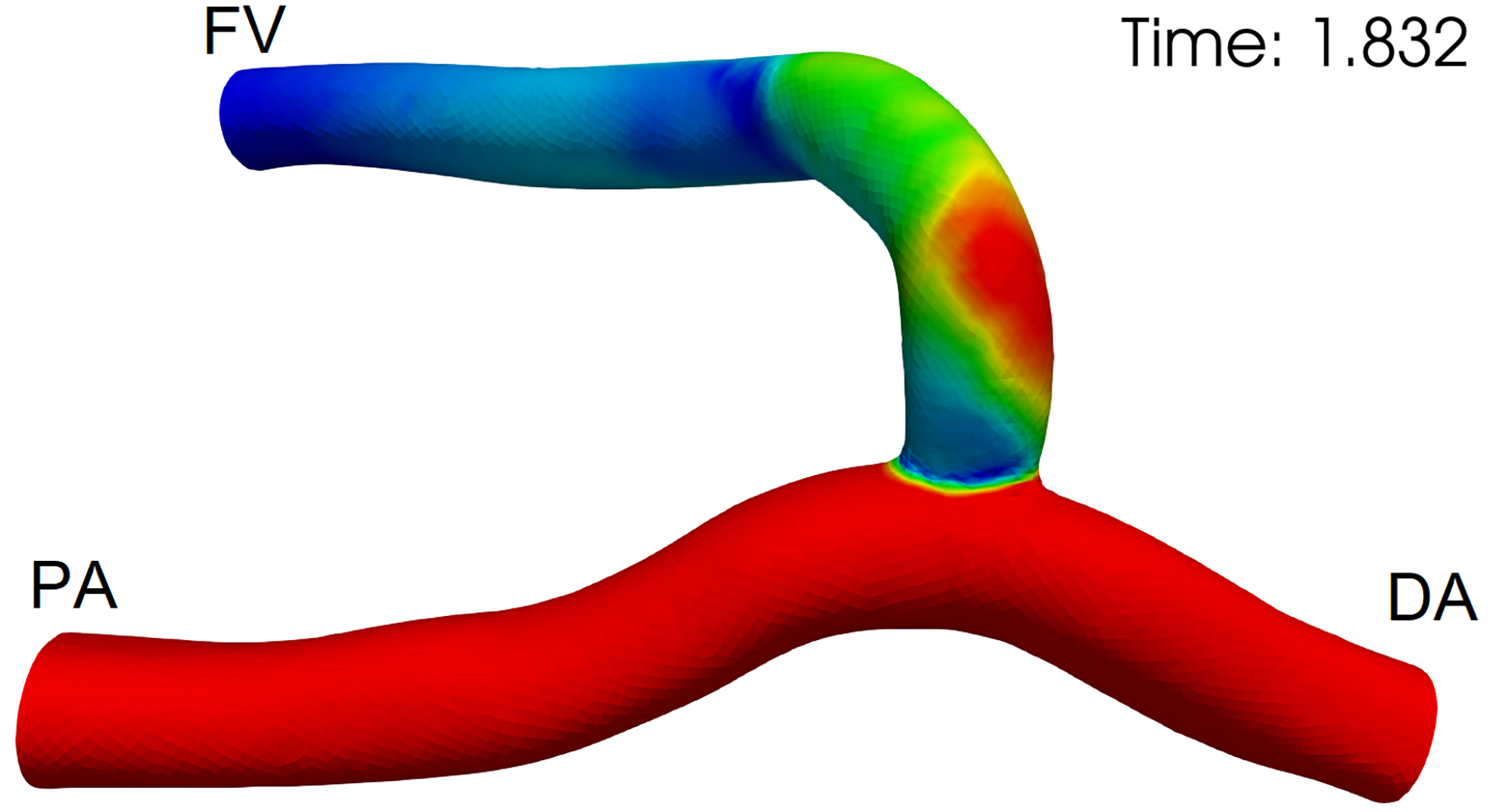}
\end{center}
\hspace{6.75cm} (c) Case O (One-way) 
\begin{center}
\includegraphics[keepaspectratio,width=7.75cm]{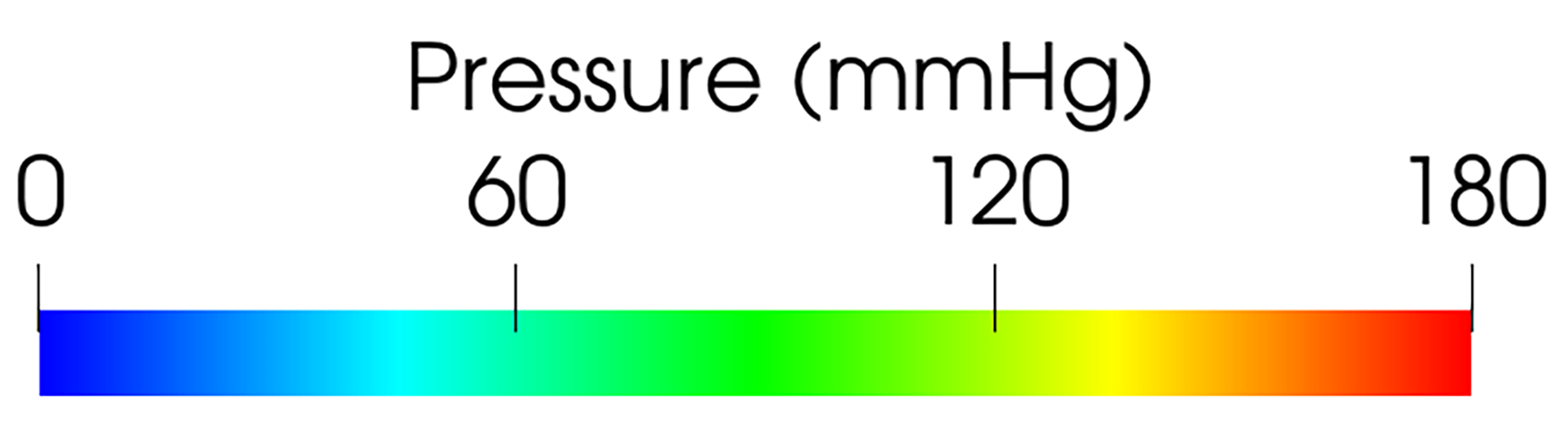}
\end{center}
\caption{Pressure differences for the S, M, and O cases.}
\label{fig6:press}
\end{figure}

Figure~\ref{fig6:press} presents pressure difference distributions, where the reference pressure value is set at zero on the outflow boundary; the pressure difference is higher in cases M and O than in case S. In case O depicted in Fig.~\ref{fig6:press}(c), a higher-pressure region arises near the anastomotic site in FV, which is attributable to the circulation induced by a small disturbance through regurgitation from DA.

\subsection{Energy loss and pressure drops}
\label{sec3.3:Energy}
\begin{figure}[htbp]
\begin{center}
\includegraphics[keepaspectratio,width=8.cm]{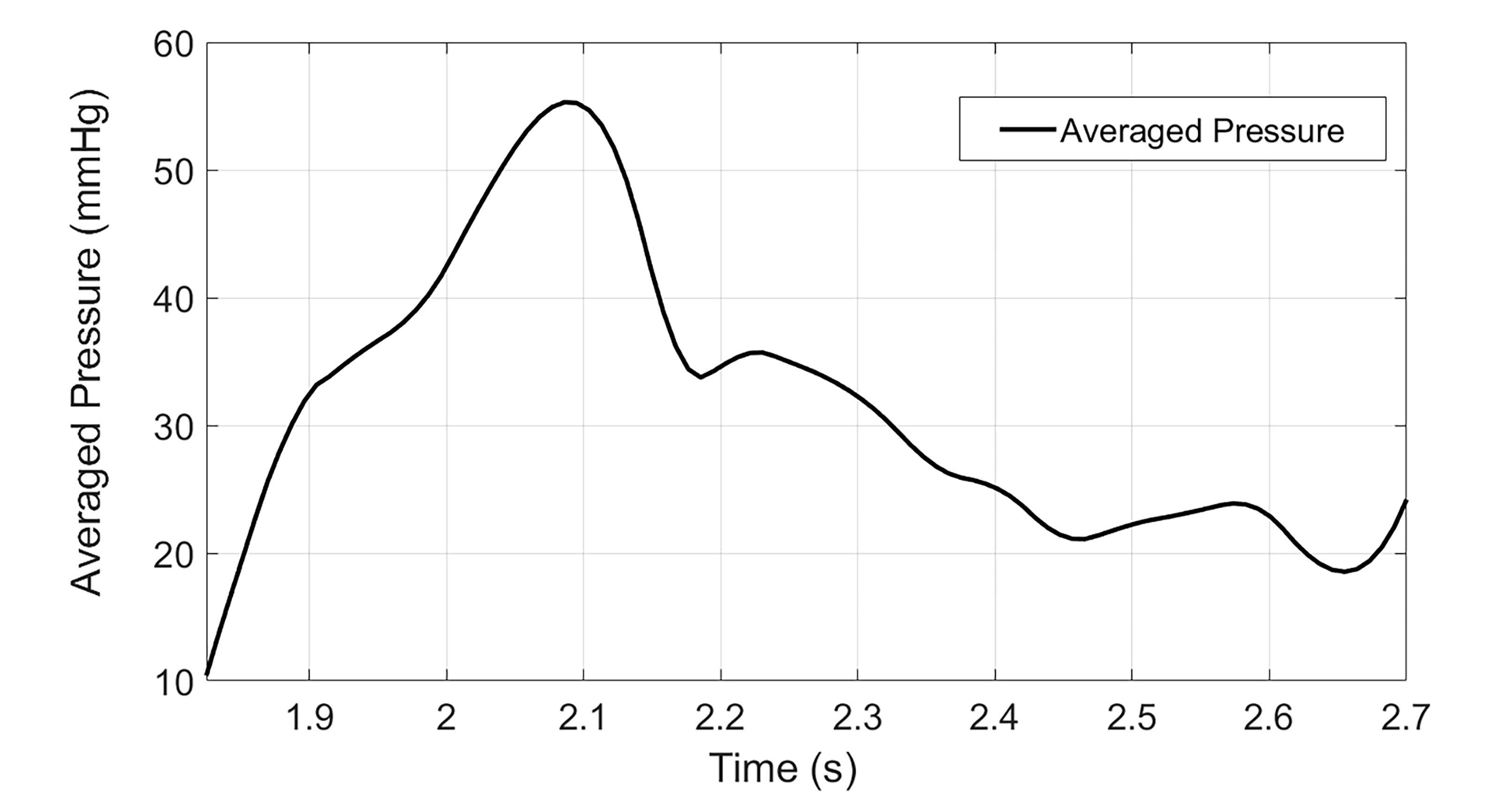}
\includegraphics[keepaspectratio,width=8.25cm]{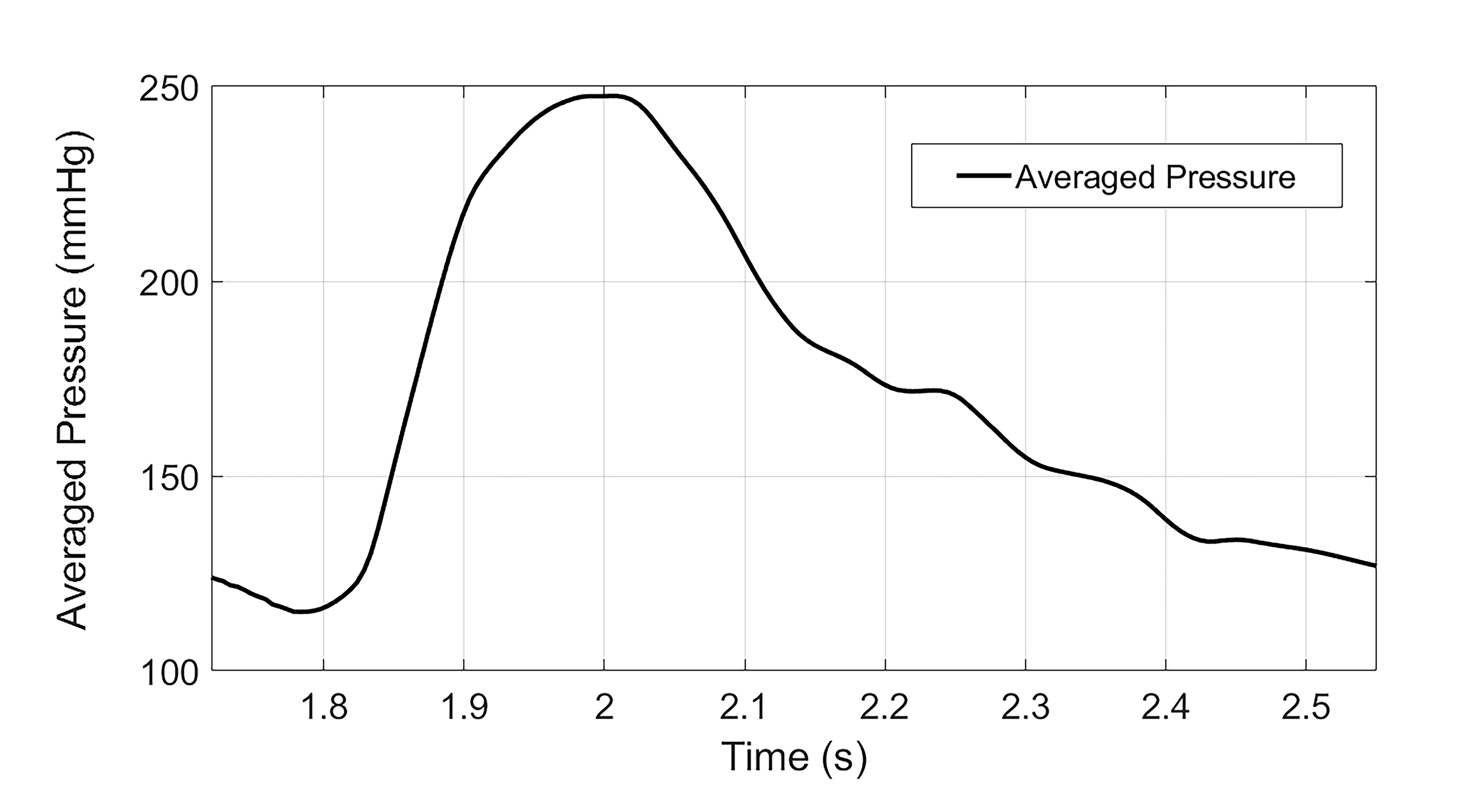}
\end{center}
\hspace{2.5cm} (a) Case S (Splitting)  \hspace{5.25cm} (b) Case M (Merging)
\begin{center}
\includegraphics[keepaspectratio,width=8.25cm]{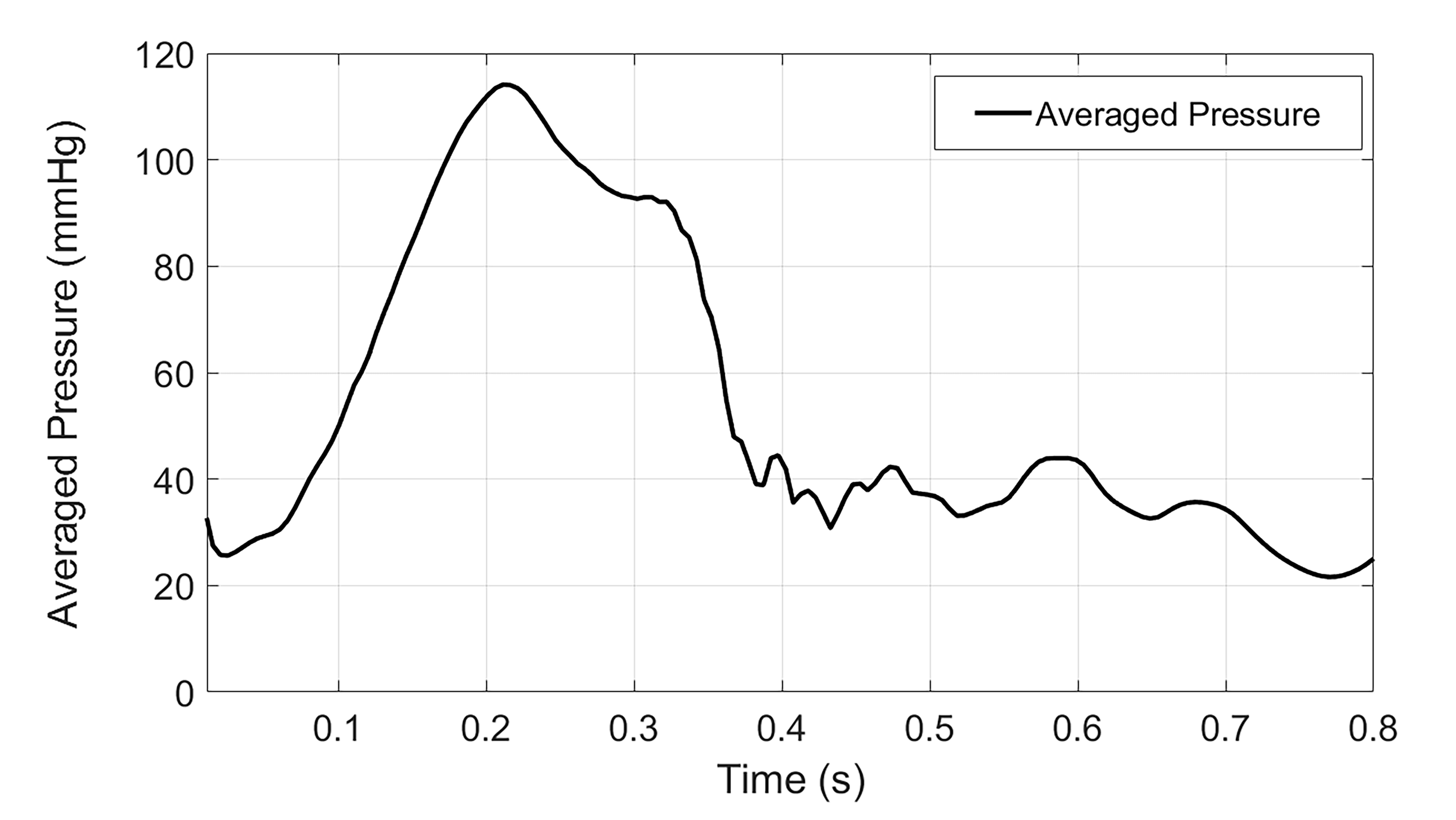}
\end{center}
\hspace{7cm} (c) Case O (One-way) 
\caption{Spatially Averaged Pressures for the S, M, and O cases.}
\label{fig7:avgpress}
\end{figure}

Figure~\ref{fig7:avgpress} presents the spatially averaged pressure history for one heart period in cases S, M, and O, respectively, which shows that the pressures change greatly in cases M and O compared with case S. Differences of time-averaged pressures between inflow and outflow boundaries are 45 mmHg, 130 mmHg, and 90 mmHg, respectively, in cases S, M, and O. The pressure drops between inflow and outflow boundaries are related strongly to the total flow rate as well as disturbances in the flow fields near the anastomotic site. The disturbed flow in FV near the anastomotic site in case O seems to contribute to larger pressure drop in this case. From the perspective of energy balance, a larger pressure drop signifies a stronger energy loss, which affects degeneration of endothelial cells of vessel walls through the effects of WSS \cite{browne2015role,ene2012disturbed,ku1985pulsatile}.

\section{Discussion}
\label{sec4:Dis}
Patients with ESRD often need hemodialysis with appropriate VAs, of which the most preferred type is AVFs. Prolonged patency of the AVFs is strongly anticipated but approximately 50\% of AVFs fail sometimes after surgery, causing adverse clinical outcomes \cite{friedman1986shear,corpataux2002low}. The underlying failure mechanism remains to be elucidated, but the most common reason for failure is the formation of intimal hyperplasia, which reduces the size of the lumen and flow rate within AVF. For this study, we have considered three patient cases including the regurgitating flow case, in which the blood flow regurgitates to the radial artery from the ulnar artery through the palmar arch because of the changes in pressure balance between radial and ulnar arteries. Although we were unable to measure the actual flow rate in the distal part (DA) correctly because of thrill of the vessels, we have applied the mass conservation law to obtain more accurate flow rate data.

In case M (merging), the regurgitating flow strongly influences the main flow at the anastomotic site. The strong flow adheres progressively to the outer wall of FV, elevating the shear stress there as presented in Fig.~\ref{fig5:cross}(b). In case O (one-way), a small flow from DA disturbs the main flow in FV. Circulations around the vessel axis are apparent in FV near the anastomotic site shown in Fig.~\ref{fig3:streamlines}(c). These results exhibited that the M case with strong regurgitation is subjected to quite high shear stresses, which concurs with results reported from several earlier studies \cite{ene2001,kharboutly2007,kharboutly2010numerical,huijbregts2007hospital,ene2012disturbed}, and indicating that vessel walls are more vulnerable to vascular damage when the WSS is high or oscillating. Disturbances of flow at venous anastomosis of AVFs are also associated with vascular remodeling described in these reports \cite{fillinger1989beneficial,ku1985pulsatile}. It has been observed that circulation in the venous side at anastomotic site, caused by regurgitating flow, is also responsible for significant pressure drops between the inflow and outflow boundaries for M and O cases. We have observed that venous segments of these AVFs showed much lower blood pressure levels than the arterial segments depicted in Fig.~\ref{fig6:press}. These computational results are in line with those of earlier studies \cite{browne2015vivo,stella2019assessing,bozzetto2016transitional}.

AVF patients undergoing hemodialysis usually have complicated hemodynamic conditions. The surgical construction of the anastomosis between an artery and a vein affects the flow dynamics there because of local morphologies of the blood vessels and because of global pressure distribution changes. This paper has presented computational analysis using patient-specific morphologies and flow conditions immediately after construction of AVFs from fluid dynamics points of view. As time goes on after surgery, deformation and pressure balance in the vessels network might change, which finally brings about occlusion of vessels. The present study has the following limitations. We considered a small number of patient cases that represent flow conditions of specific types. A broad scope of AVF morphologies and regurgitating patterns must be considered in future works. Moreover, the morphologies adopted for this study were taken only a week after the surgery. Morphologies and flow conditions might change after a long time, which is beyond the scope of the present study. 

The results presented herein suggest important differences between different morphologies and flow conditions including regurgitations. Regurgitation brings about high wall shear stress near the anastomotic site because of instabilities induced by merging phenomena, for which type careful follow-up examinations are regarded as necessary. These findings might therefore provide insights into reasons for occlusions and maturation failures, along with indications for surgery planning.

\section*{Conflict of interest} We hereby declare that the authors have no potential conflict of interest related to this study. 

\section*{Ethical Statement}
Institutional review board (IRB) approval was obtained from Akashi Medical Center Hospital and the Tohoku University School of Medicine in accordance with the Declaration of Helsinki. 

\section*{Funding} 
This work was supported by JST CREST with Grant Number: JPMJCR15D1, Japan.

\singlespacing
\bibliographystyle{unsrturl}  
\bibliography{references}

\end{document}